\documentclass[letterpaper, 12pt]{amsart}

\usepackage[all]{xy}                        
\UseTips                                    

\usepackage{amssymb}                        
\usepackage{xspace}                         
\usepackage{stmaryrd}                       
\usepackage[margin=1.2in]{geometry}         

\theoremstyle{plain}
\newtheorem{theorem}{Theorem}[section]
\newtheorem*{theorem*}{Theorem}
\newtheorem{proposition}[theorem]{Proposition}
\newtheorem{corollary}[theorem]{Corollary}
\newtheorem{lemma}[theorem]{Lemma}
\newtheorem{conjecture}[theorem]{Conjecture}
\newtheorem{exercise}{Exercise}[section]

\theoremstyle{definition}
\newtheorem{definition}[theorem]{Definition}

\theoremstyle{remark}
\newtheorem{remark}[theorem]{Remark}
\newtheorem{example}[theorem]{Example}

\newcommand{\enm}[1]{\ensuremath{#1}}          
\newcommand{\op}[1]{\operatorname{#1}}
\newcommand{\cal}[1]{\mathcal{#1}}

\newcommand{\floor}[1]{\llcorner #1\lrcorner}

\newcommand{\CC}{\enm{\mathbb{C}}}             
\newcommand{\NN}{\enm{\mathbb{N}}}
\newcommand{\RR}{\enm{\mathbb{R}}}
\newcommand{\QQ}{\enm{\mathbb{Q}}}
\newcommand{\ZZ}{\enm{\mathbb{Z}}}

\renewcommand{\AA}{\enm{\mathbb{A}}}
\newcommand{\PP}{\enm{\mathbb{P}}}
\newcommand{\LL}{\enm{\mathbb{L}}}

\newcommand{\Ff}{\enm{\cal{F}}}

\newcommand{\Ii}{\enm{\cal{I}}}
\newcommand{\Jj}{\enm{\cal{J}}}

\newcommand{\Mm}{\enm{\cal{M}}}

\newcommand{\Oo}{\enm{\cal{O}}}

\renewcommand{\phi}{\varphi}        
\renewcommand{\theta}{\vartheta}
\renewcommand{\epsilon}{\varepsilon}

\newcommand{\Spec}{\op{Spec}}

\newcommand{\Hom}{\op{Hom}}

\newcommand{\codim}{\op{codim}}

\newcommand{\invlim}{\varprojlim}

\newcommand{\Image}{\op{Im}}
\newcommand{\ord}{\op{ord}}

\newcommand{\tensor}{\otimes}         

\renewcommand{\to}[1][]{\xrightarrow{\ #1\ }}

\newcommand{\into}{\hookrightarrow}

\newcommand{\diff}[1][x]{\textstyle{\frac{\partial}{\partial #1}}}

\newcommand{\usc}[1][m]{\underline{\phantom{#1}}}

\newcommand{\xfootnote}[1]{}
\newcommand{\solution}[1]{\footnote{#1}}

\newcommand{\defeq}{\stackrel{\scriptscriptstyle \op{def}}{=}}

\newcommand{\ie}{\textit{i.e.}\ }           
\newcommand{\eg}{\textit{e.g.}\ }           
\newcommand{\cf}{\textit{cf.}\ }

\newcommand{\jet}[2][m]{\mathcal{J}_{#1}(#2)}
\newcommand{\jetinf}[1]{\mathcal{J}_{\infty}(#1)}

\newcommand{\Mmhat}{\hat{\Mm}}
\newcommand{\Var}{\op{Var}}
\newcommand{\Der}{\op{Der}}
\newcommand{\lbr}{\llbracket}
\newcommand{\rbr}{\rrbracket}
\newcommand{\Mustata}{Musta\c{t}\v{a}\xspace}
\newcommand{\Bl}{\op{Bl}}
\newcommand{\Cont}{\op{Cont}}
\newcommand{\lct}{\op{lct}}
\newcommand{\Sing}{\op{Sing}}
\newcommand{\mld}{\op{mld}}

\begin{document}

\title{A short course on geometric motivic integration}
\author{Manuel Blickle}
\address{Universit\"at Duisburg--Essen, Standort Essen, 45117 Essen,
Germany}
\email{manuel.blickle@uni-essen.de}
\urladdr{www.mabli.org}
\date{28. July, 2005}

\begin{abstract}
These notes grew out of the authors effort to understand the theory
of \emph{motivic integration}. They give a short but thorough
introduction to the flavor of motivic integration which nowadays
goes by the name of \emph{geometric motivic integration}. Motivic
integration was introduced by Kontsevich and the foundations were
worked out by Denef, Loeser, Batyrev and Looijenga. We focus on the
smooth complex case and present the theory as self contained as
possible. As an illustration we give some applications to birational
geometry which originated in the work of \Mustata.
\end{abstract}

\maketitle

\tableofcontents

\section{The invention of motivic integration.}

 Motivic integration was
introduced by Kontsevich \cite{KonOrsay} to prove the following result
conjectured by Batyrev: Let
\[
\xymatrix{    {X_1} \ar[rd]_{\pi_1} && {X_2}\ar[ld]^{\pi_2} \\
                            &{X}&
         }
\]
be two crepant resolutions of the singularities of $X$, which itself
is a complex projective Calabi-Yau\footnote{Usually, a normal
projective variety $X$ of dimension $n$ is called Calabi-Yau if the
canonical divisor $K_X$ is trivial and $H^i(X,\Oo_X)=0$ for $0<i<n$.
This last condition on the cohomology vanishing is not necessary for
the statements below. In the context of mirror symmetry it seems
customary to drop this last condition and call $X$ Calabi-Yau as
soon as $K_X=0$ (and the singularities are mild), see
\cite{Bat.StringyHodge}.} variety with at worst canonical Gorenstein
singularities. Crepant (as in \emph{non discrepant}) means that the
pullback of the canonical divisor class on $X$ is the canonical
divisor class on $X_i$, \ie the discrepancy divisor $E_i = K_{X_i} -
\pi_i^*K_{X}$ is numerically equivalent to zero. In this situation
Batyrev showed, using $p$-adic integration, that $X_1$ and $X_2$
have the same betti numbers $h_i=\dim H^i(\usc,\CC)$. This lead
Kontsevich to invent \emph{motivic integration} to show that $X_1$
and $X_2$ even have the same Hodge numbers $h^{i,j} = \dim
H^i(\usc,\Omega^j)$.

This problem was motivated by the \emph{topological mirror symmetry test} of
string theory which asserts that if $X$ and $X^*$ are a mirror pair\footnote{To
explain what a mirror pair is in a useful manner lies beyond my abilities. For
our purpose one can think of a mirror pair (somewhat tautologically) as being a
pair that passes the topological mirror symmetry test. Another achievement of
Batyrev \cite{Bat.MirorDual} was to explicitly construct the mirror to a mildly
singular (toric) Calabi-Yau variety.} of smooth Calabi-Yau varieties then they
have mirrored Hodge numbers
\[
    h^{i,j}(X) = h^{n-i,j}(X^*).
\]
As the mirror of a smooth Calabi-Yau might be singular, one cannot restrict to
the smooth case and the equality of Hodge numbers actually fails in this case.
Therefore Batyrev suggested, inspired by string theory, that one should look
instead at the Hodge numbers of a crepant resolution, if such exists\footnote{
Calabi-Yau varieties do not always have crepant resolutions. I think one of
Batyrev's papers discusses this.}. The independence of these numbers from the
chosen crepant resolution is Kontsevich's result. This makes the \emph{stringy
Hodge numbers} $h^{i,j}_{st}(X)$ of $X$, defined as $h^{i,j}(X')$ for a crepant
resolution $X'$ of $X$, well defined. This leads to a modified mirror symmetry
conjecture, asserting that the stringy Hodge numbers of a mirror pair are equal
\cite{Bat.MirorDual}.

Batyrev's conjecture is now Kontsevich's theorem and the simplest
form to phrase it might be:
\begin{theorem}[Kontsevich]
Birationally equivalent smooth Calabi-Yau varieties have the same
Hodge numbers.\footnote{There is now a proof by Ito
\cite{Ito.StringyHodge} of this result using $p$-adic integration,
thus continuing the ideas of Batyrev who proved the result for Betti
numbers using this technique. Furthermore the recent weak
factorization theorem of W{\l}odarczyk
\cite{Wlodar.AbretalFactorize} allows for a proof avoiding
integration of any sort.}
\end{theorem}
\begin{proof}
     The idea now is to assign to any variety a \emph{volume} in a suitable
     ring $\Mmhat_k$ such that the information about the Hodge numbers is retained. The
     following diagram illustrates the construction of $\Mmhat_k$:
     \[
     \xymatrix{
     {\Var_k} \ar[r]\ar[rd]_{E} &{K_0(\Var_k)} \ar[r]\ar[d] &{K_0(\Var_k)[\LL^{-1}]}\ar[r]\ar[d]&{\Mmhat_k} \ar[d]\\
     {}&{\ZZ[u,v]}\ar@{^(->}[r] &{\ZZ[u,v,(uv)^{-1}]}\ar@{^(->}[r] &{\ZZ[u,v,(uv)^{-1}]^\wedge}
     }
     \]
     The diagonal map is the (compactly supported) Hodge characteristic, which
     on a smooth projective variety $X$ is given by $E(X)=\sum (-1)^i\dim
     H^i(X,\Omega^j_X)u^iv^j$. In general it is defined via mixed Hodge
     structures\footnote{Recently, Bittner \cite{Bittner.univEul} gave an alternative construction of
     the compactly supported Hodge characteristic. She uses the weak factorization theorem of W{\l}odarczyk
     \cite{Wlodar.AbretalFactorize}
     to reduce the definition of $E$ to the case of $X$ smooth and projective,
     where it is as given above.}
     \cite{Deligne.Hodge1,Deligne.Hodge2,Deligne.Hodge3}, satisfies $E(X \times Y)= E(X)E(Y)$ for all varieties $X,Y$ and
     has the property that for $Y \subseteq X$
     a closed $k$-subvariety one has $E(X)=E(Y)+E(X-Y)$. Therefore the Hodge
     characteristic factors through the \emph{naive Grothendieck ring}
     $K_0(Var_k)$ which is the universal object with the latter property.\footnote{$K_0(\Var_k)$ is the free abelian group on the isomorphism classes
     $[X]$ of $k$-varieties subject to the relations $[X]=[X-Y]+[Y]$ for $Y$ a closed
     subvariety of $X$. The product is given by $[X][Y]=[X \times_k Y]$. The symbol $\LL$ denotes the class
     of the affine line $[\AA^1_k]$.} This explains the left triangle of the
     diagram.

     The bottom row of the diagram is the composition of a localization (inverting $uv$) and a
     completion with respect to negative degree. $\Mm_k$ is constructed analogously, by first inverting
     $\LL^{-1}=[\AA^1_k]$ (a pre-image of $uv$) and then completing appropriately (negative dimension).
     Whereas the bottom maps are injective (easy exercise), the map $K_0(\Var_k) \to \Mmhat_k$ is most
     likely not injective. The need to work with $\Mmhat_k$ instead of $K_0(\Var_k)$
     arises in the setup of the integration theory an will become clear later.\footnote{In fact, recent results of F. Loeser and R. Cluckers
     \cite{CluLoe}, and J. Sebag \cite{Sebag.motformel} indicate that the full completion may not be necessary, and all the volumes of
     measurable sets are contained in a subring of $\Mmhat_k$ that can be constructed explicitly.}

     Clearly, by construction it is now enough to show that birationally equivalent
     Calabi-Yau varieties have the same \emph{volume}, i.e.\ the same class
     in $\Mmhat_k$.
     This is achieved via the all important \emph{birational transformation rule} of motivic
     integration. Roughly it asserts that for a proper birational map $\pi: Y \to X$ the class
     $[X] \in \Mmhat_k$ is an \emph{expression} in
     $Y$ and $K_{Y/X}$ only:
     \[
        [X] = \int_{Y} \LL^{-\ord_{K_{Y/X}}}
     \]

     To finish off the proof let $X_1$ and $X_2$ be birationally equivalent Calabi-Yau varieties.
     We resolve the birational map to a Hironaka hut:
     \[
     \xymatrix{
     {}&{Y}\ar[ld]_{\pi_1}\ar[rd]^{\pi_2}&{} \\
     {X_1}\ar@{-->}[rr]&{}&{X_2}
     }
     \]

     By the Calabi-Yau assumption we have $K_{X_i} \equiv 0$ and therefore
     $K_{Y/X_i}\equiv K_Y- \pi_i^*K_{X_i} \equiv K_{Y}$. Hence the divisors $K_{Y/X_1}$
     and $K_{Y/X_2}$ are numerically equivalent. This numerical equivalence implies
     in fact an equality of divisors $K_{X/X_1}=K_{X/X_2}$ since, again by the
     Calabi-Yau assumption, $\dim H^0(X,K_Y)=\dim
     H^0(X_i,\Oo_{X_i})=1$.\footnote{In general, the condition that $X_1$ and $X_2$ have a
     common resolution $Y$ such that $K_{Y/X_1}$ is numerically equivalent $K_{Y/X_1}$ is
     called $K$--equivalence. We showed above that two birational Calabi--Yau varieties are $K$--equivalent.
     For mildly singular $X_i$ (say canonical) it can be derived from the negativity lemma
     \cite[Lemma 3.39]{KolMor} that $K$--equivalence implies actual equality of divisors $K_{Y/X_1}$ and $K_{Y/X_1}$.
     Hence the Calabi--Yau assumption was not essential to conclude this (but provides a simple argument).}
     By the transformation rule, $[X_1]$ is an expression depending
     only on $Y$ and $K_{X/X_1}=K_{X/X_2}$. The same is true for $[X_2]$ and
     thus we have $[X_1]=[X_2]$ as desired.
\end{proof}

These notes were started during a working seminar at MSRI during the
year of 2003 and took shape in the course of the past 2 years while
I was giving introductory lectures on the subject. They have taken
me way too much time to finish and would not have been finished at
all if it weren't for the encouragement of many people: Thanks goes
to all the participants of the seminar on motivic integration at
MSRI (2002/2003), of the Schwerpunkt Junioren Tagung in Bayreuth
(2003) and the patient listeners of the mini-courses at KTH,
Stockholm (2003), the University of Helsinki (2004) and the Vigre
graduate course in Salt Lake City (2005). Special thanks goes to
Karen Smith for encouragement to start this project and to Julia
Gordon for numerous comments, suggestions and careful reading.

\section{Geometric motivic integration}
We now assume that $k$ is algebraically closed and of characteristic
zero. In fact, there is one point (see section
\ref{sec.propMotMeas}) where we will assume that $k=\CC$ in order to
avoid some technicalities which arise if the field is not
uncountable. Thus the reader may choose to replace $k$ by $\CC$
whenever it is comforting. We stress that there are significant
(though manageable) obstacles one has to overcome if one wants to
(a) work with singular spaces or (b) with varieties defined over
fields which are not uncountable or not algebraically closed. Or put
differently: The theory develops naturally (for an algebraic
geometer), and easily, in the smooth case over $\CC$, as we hope to
demonstrate below. In order to transfer this intuition to any other
situation, nontivial results and extra care is necessary.

All the results in these notes appeared in the papers of Denef and
Loeser, Batyrev and Looijenga. Our exposition is particularly
influenced by Looijenga \cite{Looijenga.motivicmeasure} and Batyrev
\cite{Bat.StringyHodge}. Also Craw \cite{Craw.IntroMotivic} was very
helpful as a first reading. We also recommend the articles of Hales
\cite{Hales.WhatMot} and Veys \cite{Veys.ArcSpMotInt}, both explain
the connection to $p$-adic integration in detail, which we do not
discuss in these notes at all. The above mentioned references are
also a great source to learn about the various different
applications the theory had to date. We will discuss none of them
except for certain applications to birational geometry.

We will now introduce the building blocks of the theory. These are:
\begin{enumerate}
\item
    The value ring of the measure: \emph{$\Mmhat_k$, a localized and completed Grothendieck ring.}
\item
    A domain of integration: \emph{$\Jj_\infty(X)$, the space of formal arcs over $X$.}
\item
    An algebra of measurable sets of $\Jj_\infty(X)$ and a measure: \emph{cylinders/stable sets and the virtual euler characteristic.}
\item
    An interesting class of measurable/integrable functions: \emph{Contact order of an arc
    along a divisor.} \item
    A change of variables formula: \emph{Kontsevich's birational transformation rule.}
\end{enumerate}
These basic ingredients appear with variations in all versions of motivic
integration (one could argue: of any theory of integration).

\subsection{The value ring of the motivic measure}\label{sec.motMeasure} Here we already gravely
depart from any previous (classical?) theory of integration since the values of
our measure do not lie in $\RR$. Instead they lie in a huge ring, constructed
from the Grothendieck ring of varieties by a process of localization and
completion. This ingenious choice is a key feature of the theory.

We start with the naive Grothendieck ring of the category of varieties over
$k$.\footnote{Alternatively, the Grothendieck ring of finite type schemes over
$k$ leads to the same object because $X-X_{\op{red}}=\emptyset$. As Bjorn
Poonen points out, the finite type assumption is crucial here. If one would
allow non finite type schemes, $K_0(\Var_k)$ would be zero. For this let $Y$ be
any $k$--scheme and let $X$ be an infinite disjoint union of copies of $Y$.
Then $[X]+[Y]=[X]$ and therefore $[Y]=0$.} This is the ring $K_0(\Var_k)$
generated by the isomorphism classes of all finite type $k$--varieties and with
relation $[X]=[Y]+[X-Y]$ for a closed $k$-subvariety $Y \subseteq X$, that is
such that the inclusion $Y \subseteq X$ is defined over $k$. The square
brackets denote the image of $X$ in $K_0(\Var/k)$. The product structure is
given by the fiber product, $[X]\cdot[Y] = [X \times_k Y](=[(X \times_k
Y)_{\op{red}}])$. The symbol $\LL$ is reserved for the class of the affine line
$[\AA^1_k]$ and $1=1_k$ denotes $\Spec{k}$. Thus, for example, $[\PP^n] =
\LL^n+\LL^{n-1}+\ldots + \LL + 1$.

Roughly speaking the map $X \mapsto [X]$ is robust with respect to chopping up
$X$ into a disjoint union of locally closed subvarieties.\footnote{This is
elegantly illustrated in the article of Hales \cite{Hales.WhatMot} which
emphasizes precisely this point of $K_0(\Var_k)$ being a \emph{scissor group}.}
By using a stratification of $X$ by smooth subvarieties, this shows that
$K_0(\Var_k)$ is generated by the classes of smooth varieties.\footnote{In
\cite{Bittner.univEul}, Bittner shows that $K_0(\Var_k)$ is the abelian group
generated by smooth projective varieties subject to a class of relations which
arise from blowing up at a smooth center: If $Z$ is a smooth subvariety of $X$,
then the relation is $[X]-[Z]=[\Bl_Z X]-[E]$, where $E$ is the exceptional
divisor of the blowup.} In a similar fashion one can assign to every
constructible subset $C$ of $X$ a class $[C]$ by expressing $C$ as a
combination of subvarieties.
\begin{exercise}
    Verify the claim in the last sentence. That is: show that the map $Y
    \mapsto [Y]$ for $Y$ a closed subvariety of $X$ naturally extends to the
    algebra of constructible subsets of $X$.
\end{exercise}
\begin{exercise}\label{ex.locfibgroth}
    Let $Y \to X$ be a piecewise trivial fibration  with constant fiber
    $Z$. This means one can write $X=\bigsqcup X_i$ as a finite disjoint union of
    locally closed subsets $X_i$ such that over each $X_i$ one has $f^{-1}X_i \cong X_i \times Z$ and
    $f$ is given by the projection onto $X_i$. Show that in $K_0(\Var_k)$ one has $[Y] = [X]\cdot[Z]$.
\end{exercise}
There is a natural notion of dimension of an element of $K_0(\Var_k)$. We say
that $\tau \in K_0(\Var_k)$ is \emph{$d$--dimensional} if there is an
expression in $K_0(\Var_k)$
\[
    \tau = \sum a_i [X_i]
\]
with $a_i \in \ZZ$ and $k$-varieties $X_i$  of dimension $\leq d$, and if
there is no expression like this with all $\dim X_i \leq d-1$. The
dimension of the class of the empty variety is set to be $-\infty$. It is easy to
verify (exercise!) that the map
\[
    \dim: K_0(\Var_k) \to \ZZ\cup \{-\infty\}
\]
satisfies $\dim (\tau \cdot \tau') \leq \dim \tau + \dim \tau'$ and $\dim
(\tau+\tau') \leq \max \{\dim \tau,\dim \tau'\}$ with equality in the latter if
$\dim \tau \neq \dim \tau'$.

The dimension can be extended to the localization $\Mm_k
\defeq K_0(\Var_k)[\LL^{-1}]$ simply by demanding that $\LL^{-1}$ has dimension
$-1$. To obtain the ring $\Mmhat_k$ in which the desired measure will take
values we further complete $\Mm_k$ with respect to the filtration induced by
the dimension.\footnote{In Looijenga \cite{Looijenga.motivicmeasure}, this is
called the virtual dimension. As described by Batyrev, composing the dimension
$\dim: \Mm_k \to \ZZ\cup\{-\infty\}$ with the exponential $\ZZ \subseteq \RR
\to[\exp(\usc)] \RR_{+}$ and by further defining $\emptyset \mapsto 0$ we get a
map
\[
    \delta_k : \Mm_k \to \RR_{+,0}
\]
which is a \emph{non-archimedian norm}. That means the following properties
hold:
\begin{enumerate}
    \item $\delta_k(A) = 0$ iff $A=0=[\emptyset]$ in $\Mm_k$.
    \item $\delta_k(A+B) \leq \max\{\, \delta_k(A), \delta_k(B)\,\}$
    \item $\delta_k(A\cdot B) \leq \delta_k(A) \cdot \delta_k(B)$
\end{enumerate}
  The ring
$\Mmhat_k$ is then the completion with respect to this norm, and therefore
$\Mmhat_k$ is complete in the sense that all Cauchy sequences uniquely converge.
The condition (2) is stronger than the one used in the definition of an
archimedian norm. This non-archimedian ingredient makes the notion of
convergence of sums conveniently simple; a sum converges if and only if the
sequence of summands converges to zero. \\ If there was an equality in
condition (3) the norm would be called \emph{multiplicative}. It is unknown
whether $\delta$ is multiplicative on $\Mm_k$. However, Poonen \cite{Poo.Groth}
shows that $K_0(\Var_k)$ contains zero divisors, thus $\delta$ restricted to
$K_0(\Var_k)$ is \emph{not} multiplicative on $K_0(\Var_k)$.} The $n$th
filtered subgroup is
\[
    \Ff^n(\Mm_k) = \{\, \tau \in \Mm_k\, |\, \dim \tau \leq
    n\,\}.
\]
This gives us the following maps which will be the basis for constructing the
sought after motivic measure:
\[
    \Var_k \to K_0(\Var_k) \to[\text{invert $\LL$}] \Mm_k \to[\usc^\wedge]
    \Mmhat_k.
\]
We will somewhat ambiguously denote the image of $X \in \Var_k$ in any of the
rings to the right by $[X]$. It is important to point out here that it is
unknown whether the completion map $\usc^\wedge$ is injective, \ie whether its
kernel, $\bigcap \Ff^d(K_0(\Var/k)[\LL^{-1}])$, is zero. It is also unknown
whether the localization is injective. \footnote{In \cite{Poo.Groth} Poonen
shows that $K_0(\Var/k)$ is not a domain in characteristic zero. It is expected
though that the localization map is not injective and that $\Mm_k$ is a domain
and that the completion map $\Mm_k \to \Mmhat_k$ is injective. But recently
Naumann \cite{Naumann.AlgInd} found in his dissertation zero-divisors in
$K_0(\Var_k)$ for $k$ a finite field and these are non-zero even after
localizing at $\LL$ -- thus for a finite field $\Mm_k$ is not a domain. For
infinite fields (\eg $k$ algebraically closed) the above questions remain
open.}

\begin{exercise}
Convergence of series in $\Mmhat_k$ is rather easy. For this observe that a
sequence of elements $\tau_i \in \Mm_k$ converges to zero in $\Mmhat_k$ if and
only if the dimensions $\dim \tau_i$ tend to $-\infty$ as $i$ approaches
$\infty$. Show that a sum $\sum_{i=0}^{\infty} \tau_i$ converges if and only
the sequence of summands converges to zero.
\end{exercise}
\begin{exercise}\label{ex.geomser}
Show that in $\Mmhat_k$ the equality $\sum_{i=0}^{\infty}\LL^{-ki} =
\frac{1}{1-\LL^{-k}}$ holds.
\end{exercise}

\subsection{The arc space $\jetinf{X}$}

Arc spaces were first studied seriously by Nash \cite{Nash.Arc} who conjectured
a tight relationship between the geometry of the arc space and the
singularities of $X$, see Ishii and Kollar \cite{IshKoll.Nash} for a recent
exposition of Nash's ideas in modern language. Recent work of \Mustata\
\cite{Mustata.JetCI} supports these predictions by showing that the arc spaces
contain information about singularities, for example rational singularities of
$X$ can be detected by the irreducibility of the jet schemes for complete
intersections. In subsequent investigations he and his collaborators show that
certain invariants of birational geometry, such as the log canonical threshold
of a pair, for example, can be read off from the dimensions of certain
components of the jet schemes, see
\cite{Mustata.SingPairJet,EinMus.LocCanThresh,EinMustYas.JetDiscr} and Section
\ref{sec.birInv} where we will discuss some of these results in detail.

Let $X$ be a (smooth) scheme of finite type over $k$ of dimension $n$. An
$m$-jet of $X$ is an order $m$ infinitesimal curve in $X$, \ie it is a morphism
\[
    \theta: \Spec k[t] /t^{m+1} \to X.
\]
The set of all $m$-jets carries the structure of a scheme
$\jet{X}$, called the $m$th \emph{jet scheme}, or space of
truncated arcs. It's characterizing property is that it is right
adjoint to the functor $\usc \times \Spec k[t]/t^{m+1}$. In other
words,
\[
    \Hom(Z \times \Spec k[t]/t^{m+1}, X) = \Hom(Z, \jet{X})
\]
for all $k$-schemes $Z$, \ie $\jet{X}$ is the scheme which
represents the contravariant functor $\Hom(\usc\times \Spec
k[t]/t^{m+1},X)$.\footnote{Representability of this functor was
proved by Greenberg \cite{Green.SchLoc1,Green.SchLoc2}; another
reference for this fact is \cite{Bosch.Neron}.} In particular this
means that the $k$--valued points of $\jet{X}$ are precisely the
$k[t]/t^{m+1}$--valued points of $X$. The so called Weil restriction
of scalars, \ie the natural map $k[t]/t^{m+1} \to k[t]/t^m$, induces
a map $\pi^m_{m-1}: \jet{X} \to \jet[m-1]{X}$ and composition gives
a map $\pi^m: \jet{X} \to \jet[0]{X}=X$. As upper indices are often
cumbersome we define $\eta_m = \pi^m$ and $\phi_m=\pi^{m}_{m-1}$.

Taking the inverse limit\footnote{For this to be defined one crucially uses
that the restriction maps are \emph{affine} morphisms.} of the resulting system
yields the definition of the \emph{infinite jet scheme}, or the \emph{arc
space}
\[
    \jet[\infty]{X} = \invlim \jet{X}.
\]
Its $k$-points are the limit of the $k$-valued points $\Hom(\Spec
k[t]/t^{m+1},X)$ of the jet spaces $\jet{X}$. Therefore they correspond to the
formal curves (or arcs) in $X$, that is to maps $\Spec k \lbr t \rbr  \to
X$.\footnote{For this observe that
$\invlim\Hom(R,k[t]/t^{m+1})\cong\Hom(R,\invlim k[t]/t^{m+1})=\Hom(R,k \lbr t
\rbr )$.} There are also maps $\pi_m: \jet[\infty]{X} \to \jet{X}$ again
induced by the truncation map $k \lbr t \rbr  \to k \lbr t \rbr /t^{m+1}$. If
there is danger of confusion we sometimes decorate the projections $\pi$ with
the space. The following picture should help to remember the notation.
\begin{equation}
\xymatrix{ {\jetinf{X}} \ar[r]^{\pi_{a}} & {\jet[a]{X}}
\ar[r]^{\pi^{a}_b} & {\jet[b]{X}}
\ar[r]_{\quad\eta_b}^{\quad\pi^b} &{X} }
\end{equation}
are the maps induced by the natural surjections
\[
\xymatrix{ {k \llbracket t \rbr} \ar[r] &{k \lbr t \rbr /t^{a+1}} \ar[r] &{k \lbr t \rbr /t^{b+1}}
\ar[r] &{k}}.
\]
\begin{example}\label{ex.JetAn}
Let $X=\Spec k[x_1,\ldots,x_n] = \AA^n$. Then, on the level of $k$-points one
has
\[
    \jet{X} = \{\, \theta: k[x_1,\ldots,x_n] \to k \lbr t \rbr /t^{m+1}
    \}
\]
Such a map $\theta$ is determined by its values on the $x_i$'s, \ie
it is determined by the coefficients of $\theta(x_i)=\sum_{j=0}^m
\theta_i^{(j)}t^j$. Conversely, any choice of coefficients
$\theta_i^{(j)}$ determines a point in $\jet{\AA^n}$. Choosing
coordinates $x_i^{(j)}$ of $\jet{X}$ with
$x_i^{(j)}(\theta)=\theta_i^{(j)}$ we see that
\[
    \jet{X} \cong \Spec
    k[x_1^{(0)},\ldots,x_n^{(0)},\ldots\ldots,x_1^{(m)},\ldots,x_n^{(m)}]
    \cong \AA^{n(m+1)}.
\]
Furthermore observe that, somewhat intuitively, the truncation map $\pi^m:
\jet{X} \to X$ is induced by the inclusion \[k[x_1,\ldots,x_n] \into
k[x_1^{(0)},\ldots,x_n^{(0)},\ldots\ldots,x_1^{(m)},\ldots,x_n^{(m)}]\] sending
$x_i$ to $x_i^{(0)}$.
\end{example}

\begin{exercise}\label{exe.EquationsJetAn}
Let $Y \subseteq \AA^n$ be a hypersurface given by the vanishing of
one equation $f=0$. Show that $\jet{Y} \subseteq \jet{\AA^n}$ is
given by the vanishing of $m+1$ equations $f^{(0)},\ldots,f^{(m)}$
in the coordinates of $\jet{\AA^n}$ described above. (Observe that
$f^{(0)}=f(x^{(0)})$ and $f^{(1)}=\sum
{\diff[x_i]{f}}(x^{(0)})x_i^{(1)}$). Show that
\begin{enumerate}
\item $\jet{Y}$ is pure dimensional if and only if $\dim \jet{Y} = (m+1)(n-1)$,
    in which case $\jet{Y}$ is a complete intersection.
\item $\jet{Y}$ is irreducible if and only if $\dim
(\pi^{m}_Y)^{-1}(Y_{\op{Sing}}) < (m+1)(n-1)$.
\end{enumerate}
Similar statements hold if $Y$ is locally a complete intersection.
\end{exercise}

The existence of the jet schemes in general (that is to show the
representability of the functor defined above) is proved, for example, in
\cite{Bosch.Neron}. From the very definition one can easily derive the
following \'etale invariance of jet schemes, which, together with the example
of $\AA^n$ above gives us a pretty good understanding of the jet schemes of a
smooth variety.

\begin{proposition}\label{prop.etale}
Let $X \to Y$ be \'etale, then $\jet{X} \cong \jet{Y} \times_Y X$.
\end{proposition}
\begin{proof}
We show the equality on the level of the corresponding functors of points
\[
    \Hom(\usc,\jet{X}) \cong \Hom(\usc \times_k \Spec
    \frac{k \lbr t \rbr }{t^{m+1}},X)
\]
and
\[
\begin{split}
    \Hom(\usc,\jet{Y} \times_Y X) &= \Hom(\usc,\jet{Y}) \times \Hom(\usc,X)
    \\
    &=\Hom(\usc \times_k \Spec
    \frac{k \lbr t \rbr }{t^{m+1}},Y) \times \Hom(\usc,X).
\end{split}
\]
For this let $Z$ be a $k$--scheme and consider the diagram
\[
\xymatrix{
    {X} \ar[r] & {Y}  \\
    {Z\times \Spec k} \ar[u]^{\rm{p}}\ar[r]      & {Z \times \Spec
    \frac{k \lbr t \rbr }{t^{m+1}}}\ar[u]_{\theta}\ar@{-->}[ul]_{\tau}
    }
\]
to see that a $Z$-valued $m$-jet $\tau \in \Hom(Z \times_k \Spec \frac{k \lbr t
\rbr }{t^{m+1}},X)$ of $X$ induces a $Z$-valued $m$-jet $\theta \in \Hom(Z
\times_k \Spec \frac{k \lbr t \rbr }{t^{m+1}},Y)$ and a map $p \in \Hom(Z,X)$.
Virtually by definition formally \'etaleness \cite[Definition (17.1.1)]{EGA4.4}
for the map from $X$ to $Y$, the converse holds also, \ie $\theta$ and $p$
together induce a unique map $\tau$ as indicated.
\end{proof}

Using this \'etale invariance of jet schemes the computation carried out for
$\AA^n$ above holds locally on any smooth $X$. Thus we obtain:
\begin{proposition}\label{prop.locbun}
    Let $X$ be a smooth $k$-scheme of dimension $n$. Then $\jet{X}$ is locally an
    $\AA^{nm}$--bundle over $X$. In particular $\jet{X}$ is smooth
    of dimension $n(m+1)$. In the same way, $\jet[m+1]{X}$ is locally an
    $\AA^n$--bundle over $\jet{X}$.
\end{proposition}
Note that this is not true for a singular $X$ as can be seen already by looking
at the tangent bundle $TX = \jet[1]{X}$ which is well-known to be a bundle if
and only if $X$ is smooth. In fact, over a singular $X$ the jet schemes need
not even be irreducible nor reduced and can also be badly singular.

\subsection{An algebra of measurable sets}
The prototype of a measurable subset of $\jetinf{X}$ is a \emph{stable set}.
They are defined just right so that they receive a natural volume in $\Mm_k$.
\begin{definition}
A subset $A \subseteq \jetinf{X}$ is called \emph{stable} if for all $m \gg 0$,
$A_m
\defeq \pi_m(A)$ is a constructible subset\footnote{The constructible
subsets of a scheme $Y$ are the smallest algebra of sets containing the closed
sets in Zariski topology.} of $\Jj_m(X)$, $A=\pi_m^{-1}(A_m)$ and the map
\begin{equation}\label{eqn.defstab}
    \pi^{m+1}_m: A_{m+1} \to A_m\quad \text{is a locally trivial
    $\AA^n$--bundle.}
\end{equation}
For any $m \gg 0$ we define the \emph{volume} of the stable set $A$ by
\[
    \mu_X(A) = [A_m]\cdot \LL^{-nm} \in \Mm_k.
\]
That this is independent of $m$ is ensured by condition
\eqref{eqn.defstab} which implies that $[A_{m+1}]=[A_m]\cdot
\LL^n$.\footnote{Reid \cite{Reid.MaKay}, Batyrev
\cite{Bat.StringyHodge} and Looijenga
\cite{Looijenga.motivicmeasure} use this definition which gives the
volume $\mu_X(\jet{X}) \in \Mm_k$ of $X$ virtual dimension $n$.
Denef, Loeser \cite{Den.GemArcSpaces} and Craw
\cite{Craw.IntroMotivic} use an additional factor $\LL^{-n}$ to give
$\mu_X(\jet{X})$ virtual dimension $0$. It seems to be essentially a
matter of taste which definition one uses. Just keep this in mind
while browsing through different sources in the literature to avoid
unnecessary confusion.}
\end{definition}
Assuming that $X$ is smooth one uses Proposition \ref{prop.locbun}
to show that the collection of stable sets forms an algebra of sets,
which means that $\jetinf{X}$ is stable and with $A$ and $A'$ stable
the sets $\jetinf{X}-A$ and $A \cap A'$ are also stable. The
smoothness of $X$ furthermore warrants that so called \emph{cylinder
sets} are stable (a cylinder is a set $A=\pi_m^{-1}B$ for some
constructible $B \subseteq \Jj_m(X)$). Thus in the smooth case
condition \eqref{eqn.defstab} is superfluous whereas in general it
is absolutely crucial. In fact, a main technical point in defining
the motivic measure on singular varieties is to show that the class
of stable sets can be enlarged to an algebra of \emph{measurable}
sets which contains the cylinders. In particular $\jetinf{X}$ is
then measurable. This is achieved as one would expect by declaring a
set measurable if it is approximated in a suitable sense by stable
sets. This is essentially carried out in
\cite{Looijenga.motivicmeasure}, though there are some inaccuracies;
but everything should be fine if one works over $\CC$ and makes some
adjustments following \cite[Appendix]{Bat.StringyHodge}.\footnote{Of
course Denef and Loeser also set up motivic integration over
singular spaces \cite{DenLo.GermArc} but their approach differs from
the one discussed here in the sense that they assign a volume to the
\emph{formula} defining a constructible set rather than to the set
of ($k$-rational) points itself. Thus they elegantly avoid any
problems which arise if $k$ is small.} To avoid these technicalities
we assume until the end of this section that $X$ is smooth over the
complex numbers $\CC$.

\subsection{The measurable function associated to a
subscheme}\label{sec.measFunc} From an algebra of measurable sets there arises
naturally a notion of measurable function. Since we did not carefully define
the measurable sets --- we merely described the prototypes --- we will for now
only discuss an important class of measurable functions.

Let $Y \subseteq X$ be a subscheme of $X$ defined by the sheaf of ideals $I_Y$.
To $Y$ one associates the function
\[
    \ord_Y : \jetinf{X} \to \NN \cup \{\infty\}
\]
sending an arc $\theta: \Oo_X \to k \lbr t \rbr $ to the order of vanishing of $\theta$
along $Y$, \ie to the supremum of all $e$ such that ideal $\theta(I_Y)$ of
$k \lbr t \rbr $ is contained in the ideal $(t^e)$. Equivalently, $\ord_Y(\theta)$ is
the supremum of all $e$ such that the map
\[
    \Oo_X \to[\theta] k \lbr t \rbr  \to k \lbr t \rbr /t^{e}
\]
sends $I_Y$ to zero. Note that this map is nothing but the truncation
$\pi_{e-1}(\theta) \in \jet[e-1]{X}$ of $\theta$.\footnote{At this point we
better set $\jet[-1]{X} = X$ and $\pi_{-1} = \pi_0=\pi$ to avoid dealing with
the case $e=0$ separately.} For a $(e-1)$-jet $\gamma\in\jet[e-1]{X}$ to send
$I_Y$ to zero means precisely that $\gamma \in \jet[e-1]{Y}$.

Thus we can rephrase this by saying that $\ord_Y(\theta)$ is the supremum of
all $e$ such that the truncation $\pi_{e-1}(\theta)$ lies in $\jet[e-1]{Y}$.
Now it is clear that
\begin{eqnarray*}
    \ord_Y(\theta)\neq 0 &\Leftrightarrow& \pi(\theta) \in Y, \\
    \ord_Y(\theta) \geq s &\Leftrightarrow& \pi_{s-1}(\theta) \in
    \jet[s-1]{Y} \text{ and } \\
    \ord_Y(\theta)=\infty &\Leftrightarrow& \theta \in
    \jetinf{Y}.
\end{eqnarray*}
The functions $\ord_Y$ just introduced are examples of measurable functions
which come up in the applications of motivic integration. For a function to be
measurable, one requires that the level sets $\ord_Y^{-1}(s)$ are measurable
sets, \ie stable sets (or at least suitably approximated by stable sets). For
this consider the set $\ord_Y^{-1}(\geq s) = \{ \theta \in \jetinf{X}\,|\,
\ord_Y(\theta) \geq s\}$ consisting of all arcs in $X$ which vanish of order
\emph{at least} $s$ along $Y$. By what we just observed $\ord_Y^{-1}(\geq s) =
\pi_{s-1}^{-1}(\jet[s-1]{Y})$ is a cylinder. Therefore, the level set
$\ord_Y^{-1}(s)$ is also a cylinder equal to
\[
    \ord_Y^{-1}(\geq s) - \ord_Y^{-1}(\geq s+1) =
\pi^{-1}_{s-1}\jet[s-1]{Y} - \pi^{-1}_{s}\jet[s]{Y}.
\]
The exception is $s=0$ in which case $\ord_Y^{-1}(\geq
0)=\pi^{-1}(X)=\jetinf{X}$ and $\ord_Y^{-1}(0)=\pi^{-1}(X)-\pi^{-1}(Y)$.

Note that the level set at infinity,
\[
    \ord_Y^{-1}(\infty) = \bigcap \ord_Y^{-1}(\geq s) = \jetinf{Y}
\]
on the other hand is \emph{not} a cylinder.\footnote{This just says
that for an arc $\theta \in \jetinf{X}$ to lie in $\jetinf{Y}$ is a
condition that cannot be checked on any truncation. Exercise!}
Still, since it is the decreasing intersection of the cylinders
$\ord^{-1}_Y(\geq s+1)=\pi^{-1}_s\jet[s]{Y}$ its alleged volume
should be obtained as the limit of the volumes of these cylinders.
The volume of $\pi^{-1}_s\jet[s]{Y}$ is
$[\jet[s]{Y}]\cdot\LL^{-ns}$. The dimension of this element of
$\Mm_k$ is $\leq \dim \jet[s]{Y} - ns$. If $Y$ is smooth of
dimension $n-c$, the we have seen that $\dim \jet[s]{Y} =
(n-c)(s+1)$. Therefore, $\jetinf{Y}$ is the intersection of cylinder
sets whose volumes have dimension $\leq (n-c)(s+1)-ns = (n-c) - cs$.
For increasing $s$ these dimensions tend to negative infinity. But
recall that in the ring $\Mmhat_k$ this is exactly the condition of
convergence to zero. Thus the only sensible assignment of a volume
to $\ord_Y^{-1}(\infty)$ is zero. This argument used that $Y$ is
smooth to show that the dimension of $\jet[m]{Y}$ grows
significantly slower than the dimension of $\jet[m]{X}$. This holds
in general for singular $Y$ and we phrase it as a proposition whose
proof however is postponed until Section
\ref{sec.BitsPieces}.\label{page.Ysmooth}
\begin{proposition}
  Let $Y \subseteq X$ be a nowhere dense subvariety of $X$, then
  $\jetinf{Y}$ is measurable and has measure
  $\mu_X(\jetinf{Y})$ equal to zero.
\end{proposition}

To make this idea into a rigorous theory one has to define a larger class of
measurable subsets of $\jetinf{X}$ and this turns out to be somewhat subtle. In
section \ref{sec.BitsPieces} we will outline this briefly -- but since, no
matter what, the measure of $\ord_Y^{-1}(\infty)$ will be zero, we will move on
at this point and start to integrate.

\subsection{Definition and computation of the motivic integral} As before let
$X$ be a smooth $k$-scheme and $Y$ a subscheme. We define the \emph{motivic
integral} of $\LL^{-\ord_Y}$ on $X$ as
\[
    \int_{\jetinf{X}} \LL^{-\ord_Y} d\mu_X = \sum_{s=0}^{\infty}
    \mu(\ord_Y^{-1}(s))\cdot \LL^{-s}.
\]
Observe that the level set at infinity is already left out from this summation
as it has measure zero.

Note that the sum on the right does converge since the virtual dimension of the
summands approaches negative infinity.\footnote{The fact that for stable sets
$A \subseteq B$ we have $\dim \mu(A) \leq \dim \mu(B)$ applied to
$\ord_Y^{-1}(s) \subseteq \jetinf{X}$ gives that the dimension of
$\mu(\ord_Y^{-1}(s))\cdot \LL^{-s}$ is less or equal to $n-s$.} The notion of
convergence in the ring $\Mmhat_k$ is such that this alone is enough to ensure
the convergence of the sum. Thus it is justified to call $\LL^{-\ord_Y}$
integrable with integral as above.

It is useful to calculate at least one example. For $Y = \emptyset$
one has $\ord_Y\equiv 0$ and thus we get
\[
    \int_{\jetinf{X}} \LL^{-\ord_Y} d\mu_X = \mu(\ord_Y^{-1}(0))=[X]
\]
where we used that $X$ is smooth. A less trivial example is $Y$ a
smooth divisor in $X$. Then $\jet[s]{Y}$ is locally a
$\AA^{(n-1)s}$--bundle over $Y$. The level set is $\ord_Y^{-1}(s) =
\pi_{s-1}^{-1}\jet[s-1]{Y} - \pi_{s}^{-1}\jet[s]{Y}$ and, using that
$[\jet[s]{Y}]=[Y]\cdot\LL^{(n-1)s}$, its measure is
\[
    [\jet[s-1]{Y}]\cdot\LL^{-n(s-1)}-[\jet[s]{Y}]\cdot\LL^{-ns}
    =[Y](\LL-1)\LL^{-s}.
\]
The integral of $\ord_Y^{-1}$ is therefore
\begin{equation}\label{eqn.basic.comp}
\begin{split}
    \int_{\jetinf{X}}\LL^{-\ord_Y} d\mu_X &=
            [X-Y]+\sum_{s=1}^\infty[Y](\LL-1)\LL^{-s}\cdot\LL^{-s} \\
        &= [X-Y]+[Y](\LL-1)\LL^{-2}\sum_{s=0}^\infty\LL^{-2s} \\
        &= [X-Y]+[Y](\LL-1)\frac{1}{\LL^2(1-\LL^{-2})} \\
        &= [X-Y]+[Y](\LL-1)(\LL^2-1)^{-1} \\
        &= [X-Y]+\frac{[Y]}{\LL+1} = [X-Y]+\frac{[Y]}{[\PP^1]}
\end{split}
\end{equation}
Note the appearance of a geometric series in line 3 which is typical for these
calculations (cf.\ Exercise \ref{ex.geomser}). In fact, the motivic volumes of
a wide class of measurable subsets (namely, the {\it semi-algebraic subsets} of
Denef and Loeser in \cite{DenLo.GermArc}) belong to the ring generated by the
image of $\Mm_k$ under the completion map, and the sums of geometric series
with denominators $\LL^j$, $j>0$ \cite[Corollary 5.2]{DenLo.geomArc}.

\begin{exercise}\solution{The answer in the first case is $[X-Y]+[Y]\frac{\LL^c-1}{\LL^{c+1}-1}$. Otherwise it
can be read off from the general formula below; e.g.\ in the second
case it is $[X-D]+\frac{[D]}{[\PP^a]}$.} Compute in a similar
fashion the motivic integrals $\LL^{-\ord_Y}$ where $Y$ is as
follows.
\begin{enumerate}
    \item[a)] $Y$ is a smooth subscheme of codimension $c$ in $X$.
    \item[b)] $Y=aD$ where $D$ is a smooth divisor and $a \in \NN$.
    \item[c)] $Y=D_1+D_2$ where the $D_i$'s are smooth and in normal crossing.
    \item[d)] $Y=a_1D_1+a_2D_2$ with $D_i$ as above and $a_i$ positive integers.
\end{enumerate}
\end{exercise}

These computations are a special case of a formula which explicitly computes
the motivic integral over $\LL^{-\ord_Y}$ where $Y$ is an effective divisor
with normal crossing support.
\begin{proposition}\label{prop.comp}
Let $Y = \sum_{i=1}^s r_i D_i$ ($r_i > 0$) be an effective divisor on $X$ with
normal crossing support and such that all $D_i$ are smooth. Then
\[
    \int_{\Jj_\infty(X)} \LL^{-\ord_Y} d\mu_X = \sum _{J \subseteq \{1,\ldots,s\}}
    [D_J^\circ] (\prod_{j\in J}\frac{\LL-1}{\LL^{r_j+1}-1})
    = \sum _{J \subseteq \{1,\ldots,s\}} \frac{[D_J^\circ]}{\prod_{j\in
    J}[\PP^{r_j}]}
\]
$D_J=\bigcap_{j\in J} D_j$ (note that $D_\emptyset = X$) and $D_J^\circ = D_J -
\bigcup_{j\not\in J} D_j$.
\end{proposition}
The proof of this is a computation entirely similar to (though
significantly more complicated than) the one carried out in
\eqref{eqn.basic.comp} above; for complete detail see either Batyrev
\cite[Theorem 6.28]{Bat.StringyHodge} or Craw \cite[Theorem
1.17]{Craw.IntroMotivic}. We suggest doing it as an exercise using
the following lemma.
\begin{lemma}\label{lem.compInt}
For $J \subseteq \{1,\ldots,s\}$ redefine $r_i = 0$ if $i \not\in J$. Then
\[
    \mu(\cap \ord^{-1}_{D_j}(r_j)) = [D^0_{J}](\LL-1)^{|J|}\LL^{-\sum r_j}
\]
where $D^0_J=\cap_{j\in J} D_j - \cup_{j \not\in J} D_{J}$ and $|J|$ denotes
the cardinality of $J$.
\end{lemma}
\begin{exercise}\solution{Hint/Solution: Clearly, for an arc $\theta$ membership in $\cap
\ord^{-1}_{D_i}(r_i)$ only depends on its truncation $\pi_t(\theta)$ as long as
$t \geq \max \{\, r_i\,\}$. With the notation of Example \ref{ex.JetAn} we
write
\[
    \pi_t(\theta(x_i)) = \sum_{j=0}^t \theta^{(j)}_i t^j
\]
such that $\theta$ is determined by the coefficients $\theta^{(j)}_i$. Now we
determine what the condition $\theta \in \cap \ord^{-1}_{D_i}(r_i)$ imposes on
the coefficients $\theta^{(j)}_i$ (here it is convenient to set $r_i=0$ for $i
> s$).
\begin{enumerate}
    \item For $j < r_i$ one has $\theta^{(j)}_i = 0$.
    \item For $j = r_i$ one has $\theta^{(j)}_i \neq 0$.
    \item For $j > r_i$ one has no condition $\theta^{(j)}_i$.
\end{enumerate}
Thus the set $\pi_t(\cap_i \ord^{-1}_{D_i}(r_i))$ is the product made up form
the factors $D^0_{\op{supp} r}$ corresponding to all possible $\theta^{(0)}_i$,
a copy of $k-\{0\}$ for each possible $\theta^{(r_i)}_i$ and $r_i > 0$ and a
copy of $k$ for each $r_i<j\leq k$. Putting this together we obtain the above
formula.}
    Show that for $t \geq \max
    \{\, r_i\,\}$ one has locally an isomorphism of the $k$-points
    \[
    \pi_t(\cap_i \ord^{-1}_{D_i}(r_i)) \cong D^0_{\op{supp} r} \times
    (k-\{0\})^{|\op{supp} r|}\times k^{nt - \sum r_i}.
    \]
    Show that this implies the lemma. To prove the preceding statement
    reduce to the case that $X = U \subseteq \AA^n$ is an open subvariety of
    $\AA^n$ and $\sum D_i$ is given by the vanishing of $x_1\cdot\ldots\cdot x_s$ where
    $x_1,\ldots,x_n$ is a local system of coordinates (use Proposition \ref{prop.etale}).
    Then finish this case using the description of the arc space of $\AA^n$ as
    given in Example \ref{ex.JetAn}.
\end{exercise}

The explicit formulas of Proposition \ref{prop.comp} and Lemma
\ref{lem.compInt} are one cornerstone underlying many applications of motivic
integration. The philosophy one employs is to encode information in a motivic
integral, then using the transformation rule of the next section the
computation of this integral can be reduced to the computation of an integral
over $\LL^{-\ord_Y}$ for $Y$ a normal crossing divisor. In this case the above
formula gives the answer. Thus we shall proceed to the all important birational
transformation rule for motivic integrals.

\section{The transformation rule}\label{sec.trans.rule} The power of the theory stems from a
formula describing how the motivic integral transforms under birational
morphisms:
\begin{theorem}\label{thm.transrule}
    Let $X' \to[f] X$ be a proper birational morphism of smooth $k$-schemes and let
    $D$ be an effective divisor on $X$, then
    \[
        \int_{\jetinf{X}}\LL^{-\ord_D} d\mu_X =
        \int_{\jetinf{X'}}\LL^{-\ord_{f^{-1}D+K_{X'/X}}} d\mu_{X'}.
    \]
    As the relative canonical sheaf $K_{X'/X}$ is defined by the Jacobian ideal of
    $f$, this should be thought of as the change of variables formula for the
    motivic integral.
\end{theorem}
As a warmup for the proof, we verify the the transformation rule
first in the special case of blowing up a smooth subvariety and
$D=\emptyset$. Let $X'=\Bl_Y X$ be the blowup of $X$ along the
smooth center $Y$ of codimension $c$ in $X$. Then by \cite[Exercise
II.8.5]{Hartshorne} the relative canonical divisor is
$K_{X'/X}=(c-1)E$, where $E$ is the exceptional divisor of the
blowup. Then, using Proposition \ref{prop.comp} in its simplest
incarnation we compute
\[
\begin{split}
    \int_{\jetinf{X'}}\LL^{-\ord_{K_{X'/X}}} d\mu_{X'} &= \int_{\jetinf{X'}}\LL^{-\ord_{(c-1)E}} d\mu_{X'} \\
        &= [X'-E] + \frac{[E]}{[\PP^c]} \\
        &= [X-Y] + [Y] = [X]\,,
\end{split}
\]
where we used that $E$ is a $\PP^c$--bundle over $Y$ (by definition
of blowup) and therefore $[E]=[Y][\PP^c]$.

\subsubsection{The induced map on the arc space}\label{sec.Finfty} The proper birational map $f$ induces a map $f_\infty =
\jetinf{f}: \jetinf{X'} \to \jetinf{X}$. The first task will be to show that
away from a set of measure zero $f_\infty$ is a bijection (of sets). Let
$\Delta \subseteq X'$ be the locus where $f$ is not an isomorphism. We show
that every arc $\gamma: \Spec k\lbr t \rbr \to X$, which does not entirely lie
in $f(\Delta)$ uniquely lifts to an arc in $X'$. For illustration consider the
diagram
\[
\xymatrix{
    {\Spec (k\lbr t\rbr_{(0)})}\ar[dd]\ar@{-->}[r]\ar[rrdd]|!{[dd];[r]}\hole |!{[ddr];[r]}\hole &{X'}\ar[dd]& {(X'-\Delta)}\ar@{_{(}->}[l]\ar[dd]^{\cong}
    \\ \\
    {\Spec (k\lbr t\rbr)} \ar[r]^-{\gamma}\ar@{..>}[ruu]
    &{X} & {(X-f(\Delta))}\ar@{_{(}->}[l]
}
\]
Observe that by assumption the generic point $\Spec k\lbr
t\rbr_{(0)}$ of $\Spec k \lbr t \rbr$ does lie in $X-f(\Delta)$ and
since $f$ is an isomorphism over $X'-\Delta$ it thus lifts to $X'$
uniquely (dashed arrow). Now the valuative criterion for properness
(see \cite[Chapter 2, Theorem 4.7]{Hartshorne}) yields the unique
existence of the dotted arrow. Thus the map $f_\infty:
(\jetinf{X'}-\jetinf{\Delta}) \to (\jetinf{X} -\jetinf{f(\Delta)})$
is a bijection of $k$-valued points. Since $\jetinf{\Delta}$ has
measure zero it can be safely ignored and we will do so in the
following.
\begin{exercise}\label{exe.f_m-surj}\solution{Solution: The valuative criterion of properness shows that
  any $\gamma \not\in \jetinf{Z}$ lies in the image (where $Z \subseteq X$ is such that
  $f$ is an isomorphism from $X'-f^{-1}(Z) \to X-Z$. But for any $\gamma_m
  \in \jet{X}$ the cylinder $\pi^{-1}_m(\gamma_m)$ cannot be contained in
  $\jetinf{Z}$ since the latter has measure zero.}
  Let $f: X' \to X$ be a proper birational map (of smooth varieties).
  Then for every $m$ the map $f:\jet{X'} \to \jet{X}$ is surjective.
\end{exercise}

\begin{proof}[Proof of transformation rule]
The level sets $C^\prime_e=\ord^{-1}_{K_{X'/X}}(e)$ partition
$\jetinf{X'}$, and cutting into even smaller pieces according to order of
contact along $f^{-1}(D)$ we define
\[
    C'_{e,k}=C'_e \cap \ord_{f^{-1}D}^{-1}(k)\, \text{ and }\,
    C_{e,k}=f_\infty(C'_{e,k})
\]
to get the following partitions (up to measure zero by \ref{sec.Finfty}) of the
arc spaces
\[
    \jetinf{X'} = \bigsqcup C'_{e,k}\, \text{ and }\, \jetinf{X} = \bigsqcup
    C_{e,k}.
\]
The essence of the proof of the transformation rule is captured by
the following two crucial facts.
\begin{enumerate}
    \item[(a)] $C_{e,k}$ are stable sets for all $e,k$.
    \item[(b)] $\mu_{X}(C_{e,k})=\mu_{X'}(C'_{e,k})\cdot \LL^{-e}$.
\end{enumerate}
Two different proofs of these facts will occupy the remainder of
this section. Using these facts, the transformation rule is a simple
calculation:
\begin{equation*}
\begin{split}
    \int_{\jetinf{X}} \LL^{-\ord_D} d\mu_X &= \sum_k \mu(\ord^{-1}_D(k)) \LL^{-k} = \sum_k \left(\sum_e \mu(C_{e,k})\right) \LL^{-k} \\
                                    &= \sum_{e,k} \mu(C'_{e,k}) \LL^{-e}
                                    \LL^{-k} \\
                                    &= \sum_t \left(\sum_{e+k=t}
                                    \mu(C'_{e,k})\right)\LL^{-t}
                                    = \sum_t
                                    \mu(\ord^{-1}_{f^{-1}D+K_{X'/X}}(t))\LL^{-t}
                                    \\
                                    &= \int_{\jetinf{X'}}
                                    \LL^{-\ord_{f^{-1}D+K_{X'/X}}} d\mu_{X'}
\end{split}
\end{equation*}
Besides facts (a) and (b) one uses that if a cylinder $B$ is written
as a disjoint union of cylinders $B_i$ then the measure $\mu(B)=\sum
\mu(B_i)$. To check that this is correct one has to use the precise
definition of the motivic measure which we have avoided until now.
See the Section \ref{sec.BitsPieces} for more details on this.
\end{proof}

We first treat properties (a) and (b) in a special case, namely the
blowup at a smooth center. In the applications of motivic
integration to birational geometry which we discuss below, we can
always put ourselves in the favorable situation that the birational
map in consideration is a sequence of blowing ups along smooth
centers, hence already this simple version goes very far.
Furthermore, using the Weak Factorization Theorem of
\cite{Wlodar.AbretalFactorize} one can make a general proof of the
transformation rule by reducing to this case. However, despite the
adjective \emph{weak} in the Weak Factorization Theorem it is a very
deep and difficult result and its use in the proof of the
Transformation rule is overkill. Therefore we give in Appendix
\ref{sec.proofTrans} an essentially elementary (though at the first
reading somewhat technical) proof following
\cite{Looijenga.motivicmeasure}.

\subsection{Images of cylinders under birational maps.} We start
with some basic properties of the behaviour of cylinders under
birational morphisms, see \cite{EinLazMus.ContLoci}. These will be
useful also in Section \ref{sec.birInv}.
\begin{proposition}\label{prop.imageCyl}
    Let $f: X' \to X$ be a proper birational map of smooth
    varieties. Let $C'=(\pi^{X'}_m)^{-1}(B') \subseteq \jetinf{X'}$ be a cylinder such that\/ $B'$ is
    a union of fibers of\/ $f_m$, then\/ $C
    \defeq f_\infty(C')=(\pi^{X}_m)^{-1}(f_mB')$ is a cylinder in
    $\jetinf{X}$.
\end{proposition}
\begin{proof}
    Clearly it is enough to show that $C = \pi^{-1}_m(B)$ since
    $B=f_m(B')$ is constructible (being the image of a constructible
    set under a finite type morphism). The nontrivial inclusion is
    $\pi^{-1}_m(B) \subseteq C$. Let $\gamma \in \pi^{-1}_m(B) \subseteq
    C$ and consider for every $p \geq m$, the cylinder
    \[
        D_p \defeq (\pi^{X'}_p)^{-1}(f_p^{-1}(\pi^X_p(\gamma)))
    \]
    which is non empty by exercise \ref{exe.f_m-surj}. Clearly, $D_p \supseteq D_{p+1}$
    which implies that the $D_p$ form a decreasing sequence of nonempty cylinders. By
    Proposition \ref{prop.cylinderBaire} the intersection of all the $D_p$ is nonempty.
    Now let $\gamma' \in \bigcap D_p$ and clearly $f_\infty \gamma'
    =  \gamma$. Furthermore, since $B'$ is
    a union of fibers of $f_m$ we have $C' \supseteq D_m$ and hence
    $C' \supseteq D_p$ for all $p\geq m$. Hence $\pi_m (\gamma') \in
    B'$.
\end{proof}
\begin{exercise}\label{exe.finftysurj}\solution{Solution: The previous proposition shows that $f_\infty(\jetinf{X'})$
is a cylinder. The surjectivity of $f_m$ which was showed in
exercise \ref{exe.f_m-surj} now implies the result (in fact even the
surjectivity of $f_0=f$ is enough at this point).}
    Let $f:X' \to X$ be a proper birational map of smooth varieties.
    Then $f_\infty~:~\jetinf{X'} \to \jetinf{X}$ is surjective.
\end{exercise}

\subsubsection{The key technical result for blowup at smooth center.} We now proceed to showing
the key technical result used in the proof of the transformation
formula.

\begin{theorem}[Denef, Loeser, \cite{DenLo.GermArc}]\label{thm.maintechnical}
Let $f: X' \to X$ be a proper birational morphism of smooth
varieties. Let $C_e'= \ord_{K_{X'/X}}^{-1}(e)$ where $K_{X'/X}$ is
the relative canonical divisor and let $C_e \defeq f_\infty C_e'$.
Then, for $m \geq 2e$
\begin{enumerate}
  \item[(a')] the fiber of $f_m$ over a point $\gamma_m \in \pi_m C_e$
  lies inside a fiber of $\pi^m_{m-e}$.
  \item[(a)] $\pi_m(C_e')$ is a union of fibers of $f_m$.
  \item[(b)] $f_m: \pi_mC_e' \to \pi_m C_e$ is piecewise trivial $\AA^e$--fibration.
\end{enumerate}
\end{theorem}
\begin{corollary}
    With the notation as in the theorem, $C_e$ is a stable at level
    $m \geq 2e$ and $[C'_e]=[C_e]\LL^{e}$.
\end{corollary}
\begin{proof}[Proof of Theorem \ref{thm.maintechnical}]
  First note (a') implies (a) since $C_e'$ is stable at level $e$.
  Furthermore using the following Lemma \ref{lem.loctriv} it is
  enough to show that a fibers of $f_m$ over a point in $C_e$ is
  affine space $\AA^e$.

  We give the proof here only in the case that $X'=\Bl_YX \to X$ where
  $Y$ is a smooth subvariety of $X$, for the general argument see Appendix \ref{sec.proofTrans} below.
  Since $X$ is smooth there is an
  \'etale morphism $\phi: X \to \AA^n$ such that $\phi(Y)$ is given by
  the vanishing of the first $n-c$ coordinates on $\AA^n$ and $\phi^{-1}(\phi(Y))=Y$. By the
  \'etale invariance of jet schemes (Proposition \ref{prop.etale}) one can hence
  assume that $Y = \AA^{n-c} \subseteq \AA^n = X$.\footnote{This is not as straight froward as as it might seem.
  One has to check that $\Bl_Y X \cong X \times_\AA^n \Bl_{\phi(Y)} \AA^n$.
  Let $I$ be the ideal defining $Y$ and $\phi_*(I)$ the ideal defining $\phi(Y)$. Locally one has $\phi^*(\phi_*(I))=I$.
  The pullback of $\phi_*(I)$ to $\Bl_Y X$ is equal to the pullback of $I$ which
  is by construction of the blow-up a principal ideal sheaf. By the universal
  property of the blowup (applied to $\Bl_{\phi(Y)} \AA^n$) we get a map $\Bl_Y X \to \Bl_{\phi(Y)} \AA^n$
  and hence a map
  \[
    \Bl_Y X \to \Bl_{\phi(Y)} \AA^n \times_{\AA^n} X
  \]
  The pullback of $I$ to $\Bl_{\phi(Y)} \AA^n \times_{\AA^n} X$
  along the second projection is equal to the pullback of
  $\phi_*(I)$ along the natural map to $\AA^n$ which factors through
  the blowup. Hence this pullback is locally principal and again by the
  universal property of blowup (applied to $\Bl_Y X \to X$) we get a
  map in the opposite direction, which is easily verified to be
  inverse to the one given above. This implies in particular that $\Bl_Y
  X$ is \'etale over $\Bl_{\phi(Y)} \AA^n$ and thus
  \[
    \jet{\Bl_Y X} \cong \Bl_Y X \times_{\Bl_{\phi(Y)} \AA^n} \jet{\Bl_{\phi(Y)}
    \AA^n} \cong X \times_{\AA^n} \jet{\Bl_{\phi(Y)} \AA^n}.
  \]
  This allows one after fixing the base point $x \in X$ of $\gamma$
  to really reduce to the case of affine space.}

  To further simplify notation (and notation only) we assume that $n=3$ and
  $c=3$, that is we only have to consider the blowup of the origin in
  $\AA^3$. By definition
  \[
    X'=\Bl_0 \AA^3 \subseteq \AA^3 \times \PP^2
  \]
  is given by the vanishing of the $2 \times 2$ minors of the matrix
  \[
    \begin{pmatrix} x_0 & x_1 & x_2 \\ y_0 & y_1 & y_2 \end{pmatrix}
  \]
  where $(x_0,x_1,x_2)$ and $(z_0, z_1,z_2)$ are the coordinates on
  $\AA^3$ and the homogeneous coordinates of $\PP^2$ respectively.
  Due to the local nature of our question it is enough to consider
  on affine patch of $X'$, say the one determined by $z_0 = 1$. The equations
  of the minors then reduce to $x_1=x_0z_1$ and $x_2=x_0z_2$
  such that on this patch the map $X' \to X$ is given by the
  inclusion of polynomial rings
  \[
    k[x_0,x_0z_1,x_0z_2] \to[f] k[x_0,z_1,z_2].
  \]
  The exceptional divisor $E$ is hence given by the vanishing of
  $x_0$. The relative canonical divisor $K_{X'/X}$ is equal to $2E$
  since $\det(\op{Jac}(f))$ is easily computed to be $x_0^2$.

  Now let $\gamma' \in C_e'$. In our local coordinates, $\gamma_m'=\pi_m(\gamma')$ is
  uniquely determined by the three truncated powerseries
\begin{equation*}\label{eqn.expl}
\begin{split}
    \gamma'_m(x_0) &=t^{e/2}\sum_{i=0}^{m-e/2} a_i t^i \quad \text{with $a_0 \neq 0$}\\
    \gamma'_m(z_1) &=\sum_{i=0}^m b_i t^i \quad \text{and}\\
    \gamma'_m(z_2) &=\sum_{i=0}^m c_i t^i,
\end{split}
\end{equation*}
where the special shape of the first one comes from the condition
that $\gamma'$ has contact order with $K_{X'/X}=2E$ precisely equal
to $e$ (we only have to consider $e$ which are divisible by $c-1=2$
since otherwise $C_e'$ is empty). Its image $\gamma_m =
f_m(\gamma'_m)=\gamma'_m \circ f_m$ is analogously determined by the
three truncated powerseries
\begin{equation}\label{eqn.gamma}
\begin{split}
    \gamma_m(x_0) &=t^{e/2}\sum_{i=0}^{m-e/2} a_i t^i \quad \text{with $a_0 \neq 0$}\\
    \gamma_m(x_0z_1)&=\gamma'_m(x_0)\gamma'_m(z_1) = t^{e/2}\sum_{i=0}^{m-e/2} a_i t^i\sum_{i=0}^m b_i t^i \mod t^{m+1}\\
    \gamma_m(x_0z_2)&=\gamma'_m(x_0)\gamma'_m(z_1) = t^{e/2}\sum_{i=0}^{m-e/2} a_i t^i\sum_{i=0}^m c_i
    t^i \mod t^{m+1}
\end{split}
\end{equation}
Expanding the product of the sums in the last two equations of
\eqref{eqn.gamma} we observe that (due to the occurence of
$t^{e/2}$) the coefficients $b_{m-\frac{e}{2}+1},\ldots,b_{m}$ are
not visible in $\gamma_m(x_0z_1)$ since they only appear as
coefficients of some $t^k$ for $k>m$. Analogously,
$\gamma_m(x_0z_2)$ does not depend on
$c_{m-\frac{e}{2}+1},\ldots,c_{m}$. Conversely, given
\[
    \gamma_m(x_0z_1)= t^{e/2}\sum_{i=0}^{m-\frac{e}{2}} \beta_i t^i
\]
and knowing all $a_i$'s (which equally show up in $\gamma'_m$ and
$\gamma_m$) we can inductively recover the $b_i$'s:
\begin{equation*}
\begin{split}
    b_0&= (\beta_0)/a_0 \qquad\qquad  \text{(note: $a_0 \neq 0$)} \\
    b_1&= (\beta_1-a_1b_0)/a_0 \\
    b_2&= (\beta_2-(a_2b_0+a_1b_1))/a_0 \\
       &\vdots \\
    b_{t} &= (\beta_t-(a_tb_0+a_{t-1}b_1+\ldots+a_1b_{t-1}))/a_0
\end{split}
\end{equation*}
and this works until $t=m-\frac{e}{2}$ since $a_t$ is known for $t
\leq m-\frac{e}{2}$. The analogous statements of course hold also
for the $c_i$'s.

Summing up these observations we see that the fiber of $f_m$ over
$\gamma_m$ is an affine space of dimension $e = 2 \cdot
\frac{e}{2}$, namely it is spanned by the last $\frac{e}{2}$ of the
coefficients $b_i$ and $c_i$. This proves part (b). Furthermore, any
two $\gamma'_m$ and $\gamma''_m$ mapping via $f_m$ to $\gamma_m$
only differ in these last $\frac{e}{2}$ coefficients, hence they
become equal after further truncation to level $m-e$.\footnote{Even
to level $m-\frac{e}{2}$ in this case of blowup of a point in
$\AA^3$. In general, for a blowup of $c$ codimensional smooth center
it is truncation to level $m-\frac{e}{c-1}$ which suffices. Hence
uniformly it is truncation to level $m-e$ which works.} This shows
$(a')$ and the proof is finished.
\end{proof}
\begin{lemma}\label{lem.loctriv}
    Let $\phi:V \to W$ be a morphism of finite type schemes such
    that all fibers $\phi^{-1}(x) \cong \AA^e \times k(x)$,
    then $\phi$ is a piecewise trivial
    $\AA^e$--fibration.\footnote{In fact, it follows from Hilbert's
    Theorem 90 that a piecewise trivial $\AA^e$ fibration is
    actually locally trivial. EXPLAIN!}
\end{lemma}
\begin{proof}
We may assume that $W$ is irreducible. Then the fiber over the
generic point $\eta$ of $W$ is by assumption isomorphic to $\AA^e$.
This means that there is an open subset $U \subseteq W$ such that
$f^{-1}(U)\cong \AA^e \times U$.\footnote{The isomorphism
$\AA^e(k(\eta)) \cong \phi^{-1}(\eta) = V \times_W \Spec k(\eta)$ is
defined via some finitely many rational functions on $W$. For any $U
\subseteq W$ such that these are regular we get $\phi^{-1}(U) \cong
\AA^e \times U$.} Now $\phi$ restricted to the complement of $U$ is
a map of the same type but with smaller dimensional base and we can
finish the argument by induction. Even though we did not make this
explicit in the proof of Theorem \ref{thm.maintechnical} the
statement about the fibers being equal to $\AA^e$ holds for all
fibers (not only fibers over closed points).
\end{proof}
\begin{exercise}\solution{We may assume that $C$ is irreducible. Let $e$ be smallest such that $C \cap
\Cont^e_{K_{X'/X}}\neq 0$. This intersection is open and dense in
$C$, and hence we may replace $C$ by $C \cap \Cont^e_{K_{X'/X}}$.
Now Theorem \ref{thm.maintechnical} applies.}
    Use Proposition \ref{prop.imageCyl} and \ref{thm.maintechnical}
    to show that if $C \subseteq \jetinf{X'}$ is a cylinder, then
    the closure of $f_\infty(C) \subseteq \jetinf{X}$ is a
    cylinder, where $f : X' \to X$ is a proper birational morphism.
\end{exercise}

\subsection{Proof of transformation rule using Weak Factorization}
With the proof given so far we have the transformation formula
available for a large class of proper birational morphisms, namely
the ones which are obtained as a sequence of blowups along smooth
centers. So in particular we have the result for all resolutions of
singularities, which is the only birational morphism we will
consider in our applications later.

Let me finish by outlining how using the Weak Factorization Theorem
one can make a full proof out of this. Let us first recall the
statement:

\begin{theorem}[Weak Factorization Theorem \cite{Wlodar.AbretalFactorize}]
    Let $\phi: X' \dashrightarrow X$ be a birational map between smooth
    complete varieties over $k$ of characteristic zero. Then $\phi$
    can be factored
    
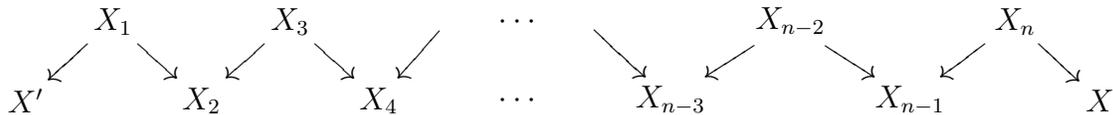
\[
    \xymatrix@R=1pc@C=1pc{
        & X_1\ar[ld]\ar[rd] && X_3\ar[ld]\ar[rd] &&\ar[dl]&\cdots&\ar[dr]&& X_{n-2}\ar[ld]\ar[rd] && X_n\ar[ld]\ar[rd] \\
        X' && X_2 && X_4&& \cdots && X_{n-3} && X_{n-1} && X
    }
    \]
    such that all indicated maps are a blowup at a smooth center.
    Furthermore, there is an index $i$ such that the rational maps
    to $X'$ to the left ($X' \dashleftarrow X_j$ for $j \leq i$) and the rational maps to $X$ to the
    right ($X_j \dashrightarrow X$ for $j \geq i$) of that index are in fact
    regular maps.
\end{theorem}

One should point out that the second part of the Theorem about the
regularity of the maps is crucial in many of its applications, and
particularly in the application that follows next.

\begin{proof}[Proof of Theorem \ref{thm.transrule}]
    So far we have proved the Transformation rule for the blowup
    along a smooth center. Given a proper birational
    morphism of smooth varieties $f: X' \to X$ we can factor it into
    a chain as in the Weak Factorization Theorem. In particular, for
    each birational map in that chain, the Transformation rule
    holds. The second part of the Factorzation Theorem together with
    the assumption that $X' \to X$ is a morphism implies that
    $X_{n-1} \to X$ is also a morphism. Part (b) of the following Exercise shows
    that for this morphism $X_{n-1} \to X$ the Transformation rule
    holds. Now the shorter chain
    ending with $X_{n-2} \to X$ is again a chain such that
    for each map the Transformation rule holds, by part (a) of that exercise. By induction we can
    conclude that the transformation rule holds for $f$ itself.
\end{proof}

\begin{exercise}\solution{Solution: Part (a) is easy. Not sure if part (c) is really feasible. Part
(d) is straightforward from (c). Let me outline part (b): Recall
that $K_{X'/X}=K_{X'/X'{}'}+{f'{}'}^*K_{X'{}'/X}$ such that
\[
\begin{split}
    \int_{X} \LL^{-\ord_{D}} d\mu_X &=
    \int_{X'}\LL^{-\ord_{{f'}^*D-K_{X'/X}}
    }d\mu_{X'} \\ &= \int_{X'}
    \LL^{-\ord_{{f'}^*D-{f'{}'}^*K_{X'{}'/X}+K_{X'/X'{}'}}} d\mu_{X'} =
    \int_{X'{}'} \LL^{-\ord_{f^*D-K_{X'{}'/X}}} d\mu_{X'{}'}
\end{split}
\]
where the first and last equality is the Transformation rule for
$f'{}'$ and $f$.} Suppose one is given a commuting diagram of proper
birational morphisms
\[
\xymatrix{
    & X' \ar[rd]^{f'}\ar[ld]_{f''}  \\
    X'' \ar[rr]^f && {X} }
\]
\begin{enumerate}
    \item[(a)] If the Transformation rule holds for $f''$ and $f$, then
    also for $f'=f'' \circ f$.
    \item[(b)] If the Transformation rule holds for $f''$ and $f'$, then it also
    holds for $f$.
    \item[(c)] The same as (a) and (b) but with ``the Transformation
    rule'' replaced by ``the conclusions of Theorem
    \ref{thm.maintechnical}''. (This part is more difficult than the
    others)
    \item[(d)] Using (c) and the Weak Factorization Theorem produce
    a proof of Theorem \ref{thm.maintechnical} building on the case
    of the blowup at a smooth center we considered above.
\end{enumerate}
\end{exercise}

\section{Brief outline of a formal setup for the motivic measure.}\label{sec.BitsPieces}
In this section we fill in some details that were brushed over in
our treatment of motivic integration so far. First this concerns
some basic properties and the well--definedness of the motivic
measure and integral.

\subsection{Properties of the motivic measure}\label{sec.propMotMeas}
For simplicity we still assume
that $X$ is a smooth $\CC$--variety. Recall that we defined for a stable set
$C=\pi_m^{-1}(B)$ a volume by setting
\[
    \mu_X(C)=[B]\LL^{-nm} \in \Mm_k.
\]
It is easy to verify that on stable sets (\ie cylinders if $X$ is smooth) the
measure is additive on finite disjoint unions. Furthermore, for stable sets $C
\subseteq C'$ one has $\dim \mu_X(C) \leq \dim \mu_X(C')$.\footnote{Check these
assertions as an exercise.} We begin with a rigorous definition of what is a
measurable set extending the above definition.
\begin{definition}\label{def.measure}
    A subset $C \subseteq \jet{X}$ is called \emph{measurable} if for all $n \in \NN$
    there is a stable set $C_n$ and stable sets $D_{n,i}$ for $i \in \NN$ such that
    \[
        C \Delta C_n \subseteq \bigcup_{i \in \NN} D_{n,i}
    \]
    and $\dim \mu(D_{n,i}) \leq -n$ for all $i$. Here $C \Delta C_n = (C-C_n) \cup
    (C_n-C)$ denotes the symmetric difference of two sets. In this case the \emph{volume} of $C$ is
    \[
        \mu_X(C) = \lim_{n \to \infty} \mu_X(C_n) \in \Mm_k.
    \]
    This limit converges and is independent of the $C_n$'s.
\end{definition}
The key point in proving the claims in the definition\footnote{In
\cite{Looijenga.motivicmeasure} a more restrictive definition of
measurable is used, namely he requires that $\dim D_{n,i} \leq
-(n+i)$. This has the advantage that one does not require the field
to be uncountable to conclude the well definedness of the measure.
Essentially, he uses that if $D \subseteq \bigcup D_i$ are cylinders
with $\lim \dim D_i = -\infty$, then $D$ is already contained the
union of finitely many of the $D_i$. This is true if $k$ is
infinite, and, if $k$ is uncountable, even true without the
assumption on the $D_i$. The advantage of Looijenga's setup is that
one is not bound to an uncountable field, but unfortunatly I was not
able to verify that in his setup $\jetinf{Y}$ is measurable and has
zero volume.} is the so called Baire property of constructible
subsets of a $\CC$--variety which crucially uses the fact that $\CC$
is uncountable, see \cite[Corolaire 7.2.6]{EGA1}.
\begin{proposition}\label{prop.baire}
Let $K_1 \supseteq K_2 \supseteq K_3 \supseteq \ldots$ be an infinite sequence
of nonempty constructible subsets of a $\CC$--variety $X$. Then $\bigcap_i K_i$
is nonempty.
\end{proposition}
For cylinder sets this implies the following Proposition, which will
be used below in a version which asserts that a cylinder $C$ which
is contained in the union of countably many cylinders $C_i$ is
already contained in the union of finitely many of these.
\begin{proposition}\label{prop.cylinderBaire}
Let $C_1 \supseteq C_2 \supseteq C_3 \supseteq \ldots$ be an
infinite sequence of nonempty cylinders in $\jetinf{X}$ where $X$ is
a smooth $\CC$--variety. Then $\bigcap_i C_i$ is nonempty.
\end{proposition}
\begin{proof}
By definition of a cylinder and Chevalley's theorem $\pi_0(C_i)$ is a
constructible subset of $X$. Thus we can apply Proposition \ref{prop.baire} to
the sequence $\pi_0(C_1) \supseteq \pi_0(C_2) \supseteq \ldots$ to obtain an
element $x_0 \in \bigcap \pi_0(C_i)$. Now consider the sequence of cylinders
$C_i'\defeq C_i \cap \pi_0^{-1}(x_0)$ and repeat the argument for the sequence
of constructible sets $\pi_1(C_1') \supseteq \pi_1(C_2')\supseteq \pi_1(C_2')
\supseteq \ldots$ to obtain an element $x_1 \in \bigcap \pi_1(C_i')$. Repeating
this procedure we successively lift $x_0\in X$ to $x_1 \in \jet[1]{X}$, $x_2
\in \jet[2]{X}$ and so forth. The limit of these $x_i$ gives an element $x \in
\jetinf{X}$ which is in the intersection of all $C_i$. Thus this intersection
is nonempty.
\end{proof}
On a possibly singular $\CC$--variety the statement (and proof) is
true if ``cylinder'' is replaced by ``stable set'' above. If $X$ is
a $k$--variety with $k$ at most countable these statements might be
false. Nevertheless, if one views the constructible sets (resp.\
cylinders) not as subsets of the $k$--rational points but rather as
certain sub--arragnements (something weaker than a subfunctor) of
the functor of points represented by $X$ the above is essentially
true. This is the point of view of Denef and Loeser and is carried
out in \cite{Den.GemArcSpaces}.

\begin{proof}[Justification of Definition \ref{def.measure}] For both claims
it suffices to show $\dim (\mu_X(C_i)-\mu_X(C'_j)) \leq -i$ for all
$j \geq i$, where the prime indicates a second set of defining data
as is in the definition. We have
\[
    C_i - C'_j \subseteq (C \Delta C_i) \cup (C \Delta C'_j) \subseteq \bigcup_m C_{i,m} \cup \bigcup_m C'_{j,m}.
\]
Since $(C_i - C'_j)$  and all terms on the right are cylinders the previous
proposition applies and $C_i - C'_j$ is contained in finitely many of the
cylinders of the right hand side. This implies that $\dim \mu_X(C_i - C'_j)
\leq -i$. The same applies to $\dim \mu_X(C'_j - C_i)$. Using that $C_i = (C_i
\cap C'_j) \cup (C_i - C'_j)$ and $C'_j = (C_i \cap C'_j) \cup (C'_j - C_i)$
and that $\mu_X$ is additive on finite disjoint unions of cylinders we get
\[
\begin{split}
    \dim (\mu_X(C_i)-\mu_X(C'_j)) &= \dim (\mu_X(C_i - C'_j) - \mu_X(C'_j - C_i)) \\
    & \leq \max \left\{ \dim \mu_X(C_i - C'_j), \dim \mu_X(C'_j - C_i) \right\}
    \leq -i.
\end{split}
\]
For once this shows that the $\mu_X(C_i)$ form a Cauchy sequence, thus the
limit exists as claimed. Secondly it immediately follows that this limit does
not depend on the chosen data.
\end{proof}
We summarize some basic properties of the measure and measurable sets.
\begin{proposition}
    The measurable sets form an algebra of sets. $\mu_X$ is additive on disjoint
    unions, thus $\mu_X$ is a pre--measure in classical terminology. If $C_i$
    are a infinite disjoint sequence of measurable sets such that
    \[
        \lim_{i \to \infty} \mu_X(C_i) = 0
    \]
    then $C= \bigcup C_i$ is measurable and $\mu_X(C)= \sum \mu_X(C_i)$.
\end{proposition}
\begin{proof}
    The verification of all parts is quite easy. As an example we only show that the
    complement of a measurable set is also measurable and leave the rest as an
    exercise.

    If $C$ is measurable then cylinders $C_i$ and
    $C_{i,j}$ can be chosen with the properties as in Definition \ref{def.measure}.
    The complements $C^c_i=\jetinf{X}-C_i$ are cylinders. Since $C^c \Delta
    C^c_j = C \Delta C_j$ it follows at once that the complement of $C$ is
    measurable.
\end{proof}
The following proposition was one of the missing ingredients for the setup of
motivic integration we outlined so far. It ensures that our typical functions
$\LL^{-\ord_Y}$ are in fact measurable (the missing part was the level set at
infinity $\jetinf{Y}$ which we owe the proof that it is measurable with measure
zero).
\begin{proposition}\label{prop.JetYzero}
    Let $Y \subseteq X$ be a locally closed subvariety. Then $\jetinf{Y}$ is a measurable
    subset of\/ $\jetinf{X}$ and if\/ $\dim Y < \dim X$ the volume
    $\mu_X(\jetinf{Y})$ is zero.
\end{proposition}
\begin{proof}\footnote{In \cite[Proposition 6.22]{Bat.StringyHodge} Batyrev claims
that with the previous results one can reduce the proof to the case
that $Y$ is a smooth divisor, where it is easily verified -- we
discussed this on page \pageref{page.Ysmooth}. Unfortunately I was
not able to follow Batyrevs argument, thus the somewhat not so self
contained proof is included here. For another another proof see
Proposition \ref{prop.dimJet}.} The proof of this result relies on a
fundamental result of Greenberg \cite{Green.RatHens} which, for our
purpose is best phrased as follows:
\begin{proposition}
    Let $Y$ be a variety. Then there exists a positive integer $c \geq 1$. such
    that
    \[
        \pi_{\floor{\textstyle\frac{m}{c}}} \jetinf{Y} = \pi^m_{\floor{\textstyle\frac{m}{c}}} \jet[m]{Y}
    \]
    for all $m \gg 0$.
\end{proposition}
This, in particular, implies that the image $\pi_n(\jetinf{Y})$ is a
constructible subset of $\jet[n]{Y}$. It can be shown that its dimension is
just the expected one, namely $\dim\pi_n (\jetinf{Y}) = (n+1)\dim Y$, see
\cite[Lemma 4.3]{Den.GemArcSpaces}. From these observations we obtain a bound
for the dimension of $\jet[m]{Y}$ as follows. We can work locally and may
assume that $Y \subseteq X$ for a smooth $X$. Then we have
\[
\begin{split}
    \dim \jet{Y} &\leq \dim (\pi^m_{\floor{\textstyle\frac{m}{c}}} \jet{Y}) + (m-\floor{\textstyle\frac{m}{c}})\dim X\\
                 &= \dim (\pi_{\floor{\textstyle\frac{m}{c}}} \jetinf{Y}) + (m-\floor{\textstyle\frac{m}{c}})\dim X
                 \\
                 &= (\floor{\textstyle\frac{m}{c}} + 1) \dim Y +(m+1)\dim X - (\floor{\textstyle\frac{m}{c}} +
                 1)\dim X \\
                 &= (m+1)\dim X - (\floor{\textstyle\frac{m}{c}}+1)(\dim X - \dim Y).
\end{split}
\]
Thus, If the codimension of $Y$ in $X$ is greater or equal to 1,
then $\dim \jet[m]{Y}\LL^{-m\dim X}$ approaches $-\infty$ as $m$
approaches $\infty$. This implies that $\jetinf{Y}$ is measurable
that its measure $\mu_X(\jetinf{Y})$ is zero.
\end{proof}

\subsubsection{Comparison with Lesbeque integration}
To guide ones intuition a comparison of the motivic measure with more classical
measures such as the Lesbeque measure or $p$--adic measures is sometimes
helpful. We discuss here the similarities with the Lesbeque measure on $\AA^n$
since this is wellknown.\footnote{The analogy with $p$--adic measures is even
more striking, see for example the discussions in \cite{Veys.ArcSpMotInt} and
\cite{Hales.PadicCompute,Hales.WhatMot}.}

For convenience we consider the case that $X=\AA^n$ in which case we identify
the $k$--points of $\jetinf{\AA_k^n}$ with $n$--tuples of power series with
coefficients in $k$. That is we identify $\jetinf{\AA_k^n} \cong (k\lbr t
\rbr)^n=\AA^n_{k\lbr t\rbr}$. Since $k\lbr t \rbr$ is a discrete valuation
domain (DVR) we can use the valuation to define a norm on $(k\lbr t \rbr)^n$ by
defining
\[
    \| \tau \| \leq \frac{1}{m} \iff  \forall i: \tau_i \in (t^m)k\lbr t \rbr
\]
for $\tau = (\tau_1,\ldots,\tau_n)$ a tuple of power series in
$\jetinf{\AA^n_k}$. It is easy to check that this defines a
(non--archimedian) norm. In Table 1 the similarities between
Lesbeques and motivic measure are summarized.

\begin{Small}
\begin{table}[h]\label{tab.lesbeques}
\begin{tabular}{|c|c|c|}
\hline
  \rule[-2mm]{0mm}{6mm}&  Lesbeques & motivic \\
\hline
 \rule[-2mm]{0mm}{6mm}space &  $\RR^n$ & $\jetinf{\AA^n}$ \\
 \rule[-2mm]{0mm}{6mm}values of measure &  $\ZZ \subseteq \QQ \subseteq \RR$ & $K_0(\Var_k) \subseteq \Mm_k \subseteq \Mmhat_k$ \\
 \rule[-2mm]{0mm}{6mm}cubes around point $a$ &  $\{x \in \RR^n | \|x-a\|\leq 1/m \}$ & $\{\gamma \in \jetinf{\AA^n}| \|\gamma-a\|\leq 1/m\}$ \\
 \rule[-2mm]{0mm}{6mm}measurable set &  $\sigma$--algebra generated by cubes & algebra of stable/measurable sets \\
 \rule[-2mm]{0mm}{6mm}volume of cube &  $(2/m)^n$ & $(\LL^{-m})^n$ \\
 \rule[-2mm]{0mm}{6mm}transformation rule & $\int_A h(f) df = \int_{g^{-1}(A)}
 h(f(x))\op{Jac}(f)dx$
  & $\int_A \LL^{-\ord_D} = \int_{f^{-1}_{\infty}(A)}
 \LL^{-\ord_D+K_{X'/X}}$ \\
\hline
\end{tabular}
\bigskip \caption{Comparison with Lesbeque measure.}
\end{table}
\end{Small}

\subsection{Motivic integration on singular varieties} In the preceding
discussion we used the assumption that your spaces are non-singular
in several places in an essential way. Rough\-ly speaking we used
the fact that if $X$ is smooth then every cylinder is a stable set,
thus can be endowed with a measure in a natural way. The fact that
the resulting algebra of measurable sets to include the cylinders
was essential for the setup since the level sets of the functions
$\ord_Y$ are cylinders in a natural way.

If $X$ is singular however, many things one might got used to from
the smooth case fail. Most prominently, the truncation maps are no
longer surjective and a cylinder is in general not stable. Thus one
has to work somewhat harder to obtain an algebra of measurable sets
which include the cylinders. In order to be able to setup an
integration theory one expects from this algebra of measurable sets
the properties asserted in the following Proposition.
\begin{proposition}
Let $X$ be a $k$--variety. Then there is an algebra of measurable subsets of
$\jetinf{X}$ and a measure $\mu_X$ on that algebra such that
\begin{enumerate}
    \item If $A$ is stable, then $A$ is measurable and $\mu_X(A)=[\pi_m A]
    \LL^{-nm}$ for $m \gg 0$.
    \item A cylinder $C$ is measurable and $\mu_X(C)= \lim_m [\pi_m C]\LL^{-nm}$.
    \item The measure is additive on finite disjoint unions.
    \item If $A \subseteq B$ are measurable, then $\dim \mu_X(A) \leq \dim
    \mu_X(B)$.
\end{enumerate}
\end{proposition}
To achieve this, one starts with the stable sets to which we know
how to assign a measure. Then one proceeds just as in the smooth
case using Definition \ref{def.measure}, and replacing ``cylinder''
by ``stable set'' whenever necessary the same proof holds as well.
The critical point now is to show that a cylinder is measurable with
volume as claimed above. Even for the cylinder $\jetinf{X} =
\pi^{-1}(X)$ it is not clear a priori that it is measurable and what
it's measure should be (in fact, the measure of $X$ to be
constructed leads to new birational invariants of $X$).

The point is that one has to partition $\jetinf{X}$ according to
intersection with the singular locus $\Sing X$, defined by the $n$th
Fitting ideal of $\Omega^1_X$. Then we can write
\[
    \jetinf{X} = \bigsqcup_{e \geq 0} \Jj^{(e)}_{\infty}(X)
\]
where $\Jj^{(e)}_{\infty}(X) = \ord^{-1}_{\op{Sing} X}(e)$. As it
turns out, the sets $\Jj^{(e)}_{\infty}(X)$ are in fact stable at
level $\geq e$. The method is analogous to partitioning the
cylinders of $X'$ according to intersection with $K_{X'/X}$ (defined
by the $0$th Fitting ideal of $\omega_{X'/X}$) in the proof of the
transformation rule. If one treats everything subordinate to this
partition according to intersection with the singular locus one can
construct a working theory in the singular case, see
\cite{Looijenga.motivicmeasure,DenLo.GermArc}.

\section{Birational invariants via motivic integration}\label{sec.birInv}
As an illustration of the theory we discuss some applications of
geometric motivic integration to birational geometry, namely we give
a description of the log canonical threshold of a pair $(X,Y)$,
where $Y$ is a closed subscheme of the smooth scheme $X$, in terms
of the asymptotic behavior of the dimensions of the jet schemes
$\jet{Y}$. These results are due to \Mustata\
\cite{Mustata.SingPairJet} and his collaborators
\cite{EinMus.LocCanThresh,EinMustYas.JetDiscr,Yasuda.DimJet}. The
precise statement is as follows:
\begin{theorem}\label{thm.lcformula}
    Let $Y \subseteq X$ be a subscheme of the smooth variety $X$. Then the log
    canonical threshold of the pair $(X,Y)$ is
    \[
        c(X,Y)=\dim X - \sup_m\left\{ \frac{\dim \jet{Y}}{m+1}\right\}.
    \]
\end{theorem}
The definition of the log canonical threshold requires the
introduction of some more notation from birational geometry which
will be done shortly. It is an invariant which can be read of from
the data of a log resolution of the pair $(X,Y)$. The proof of the
above is a very typical application of motivic integration as its
strategy is to express the quantity one is interested in (say $\dim
\jet{Y}$), in terms of a motivic integral. Then one uses the
transformation rule (Theorem \ref{thm.transrule}) to reduce to an
integral over a normal crossing divisor which can be explicitly
computed, similarly as Formula \ref{prop.comp}.

This result is in line with the earliest investigations of jet spaces by Nash
\cite{Nash.Arc} who conjectured an intimate correspondence between the jet
spaces of a singular space and the divisors appearing in a resolution of
singularities. Even though his conjecture was disproved recently by Kollar and
Ishii \cite{IshKoll.Nash} in general, there are important cases where his
prediction was true, for example in the case of toric varieties.

\subsection{Notation from birational geometry}\label{sec.notbir} Throughout this section we fix
the following setup. Let $X$ be a smooth $k$--variety of dimension $n$ and let
$k$ be of characteristic zero. Let $Y \subseteq X$ be a closed subscheme. Let
$f:X' \to X$ be a log resolution of the pair $(X,Y)$, that is a proper
birational map such that
\begin{enumerate}
    \item $X'$ is smooth and
    \item denoting by $F \defeq f^{-1}Y = \sum_{i=1}^s a_iD$ and $K \defeq K_{X'/X}
    = \sum_{i=1}^s b_iD$ for some $a_i,b_i \in \QQ$ and prime divisors $D_i$, the divisors
    $F$, $K$ and $F+K$ have simple normal crossing support.
\end{enumerate}
The existence of log resolutions is a consequence of Hironaka's
resolution of singularities. Now one defines:
\begin{definition} Let $X$ and $Y$ and the datum of a log resolution be as
above and let $q \geq 0$ be a rational number. Then we say that
\begin{enumerate}
    \item $(X, q\cdot Y)$ is \emph{Kawamata log terminal} (KLT) if and only if
    $b_i - qa_i +1 > 0$ for all $i$.
    \item $(X, q\cdot Y)$ is \emph{log canonical} (LC) if and only if
    $b_i -qa_i+1 \geq 0$ for all $i$.
\end{enumerate}
\end{definition}
We point out (without proof) that these notions are independent of
the chosen log resolution and therefore well defined, see
\cite{Laz.PAG2} for details.
\begin{remark}
These notions can be expressed in terms of the multiplier ideal $\Ii(q \cdot
I_Y)$ of $I_Y$, the sheaf of ideals which cuts out $Y$ on $X$, as follows:
\begin{eqnarray*}
(X,q\cdot Y) \text{ is KLT }\iff& \Ii(I_Y^q) = \Oo_X \\
(X,q\cdot Y) \text{ is LC\  }\iff& \Ii(I_Y^q) = \Oo_X\quad \forall
q'<q.
\end{eqnarray*}
To see this observe that by definition
\[
    \Ii(I_Y^q) = f_*\Oo_{X'}(\ulcorner K-qF \urcorner) = f_*\Oo_{X'}(\ulcorner b_i - q a_i\urcorner D_i)
\]
which is equal to $\Oo_X$ if and only if $\ulcorner b_i - q a_i\urcorner \geq
0$ for all $i$. as the upper corners denote the round up of an integer, this is
equivalent to $b_i - q a_i +1 > 0$ for all $i$ as required.
\end{remark}
Now we proceed to the definition of the log canonical threshold, which is just
the largest $q$ such that the pair $(X,q\cdot Y)$ is Kawamata log terminal.
\begin{definition}
    The \emph{log canonical threshold} of the pair $(X,Y)$ is
    \[
    \begin{split}
        \lct(X,Y) &= \sup \{\, q \,|\, (X,q \cdot Y) \text{ is KLT }\,\} \\
               &= \sup \{\, q \,|\, b_i-qa_i+1 > 0 \quad \forall i\,\} \\
               &= \min_i \left\{\, \frac{b_i+1}{a_i}\,\right\}
    \end{split}
    \]
\end{definition}
Note that clearly one has $\lct(X,q\cdot Y) = q^{-1}\lct(X,Y)$ so
that we can restrict to the case $q=1$ in the definition of the log
canonical threshold. The formula for the log canonical threshold in
terms of the jet schemes which we are aiming to proof in this
section is an immediate consequence of the following theorem.
\begin{theorem}\label{thm.KLTLC}
Let $X$ and $Y$ be as before. Then
\begin{eqnarray*}
    (X,q\cdot Y)\text{ is KLT }&\iff& \dim \jet{Y} < (m+1)(n-q)\text{ for all $m$} \\
    (X,q\cdot Y)\text{ is LC  }&\iff& \dim \jet{Y} \leq (m+1)(n-q)\text{ for all
    $m$}.
\end{eqnarray*}
\end{theorem}\label{thm.KLTChar}
From this the proof of Theorem \ref{thm.lcformula} follows
immediately.
\begin{proof}[Proof of Theorem \ref{thm.lcformula}] By definition we have
\[
\begin{split}
     \lct(X,Y) &= \sup \{\, q \,|\, (X,q \cdot Y)\text{ is KLT }\,\} \\
            &= \sup \{\, q \,|\, \dim \jet{Y} < (m+1)(n-q)\, \forall m\,\} \\
            &= n- \sup_m \left\{ \frac{\dim \jet{Y}}{m+1}\right\}
\end{split}
\]
\end{proof}
Moreover, the proof of Theorem \ref{thm.KLTLC} will reveal that the
supremum in Theorem \ref{thm.KLTChar} is actually obtained by
infinitely many $m$, namely whenever $m+1$ is divisible by all the
coefficients $a_i$ of $D_i$ in $f^{-1}Y=\sum a_iD_i$ one has
$\lct(X,Y)=\dim X -\frac{\dim \jet{Y}}{m+1}$.

Before proceeding to the proof we derive some elementary properties of the log
canonical threshold.
\begin{proposition} Let $Y \subseteq Y' \subseteq X$ closed subvarieties of
$X$. Then
\begin{enumerate}
    \item $\lct(X,Y) \geq \lct(X,Y')$.
    \item $0 < \lct(X,Y) \leq \codim(Y,X)$ with equality if $(X,Y)$ is log canonical.
    \item $\lct(X,Y)$ is independent of $X$ of fixed dimension.
    \item Let $(X',Y')$ be another pair, then $\lct(X \times X',Y \times
    Y')=\lct(X,Y)+\lct(X',Y')$.
\end{enumerate}
\end{proposition}
\begin{proof}
For (1) note that $Y \subseteq Y'$ implies that $\dim \jet{Y} \leq \dim
\jet{Y'}$ then apply Theorem \ref{thm.KLTChar}.

For (2) recall that in any case $\jet{Y}$ contains $\jet{Y_\text{reg}}$ the jet
scheme over the regular locus of $Y$. The latter has dimension $(m+1)\dim Y$.
Therefore $\dim \jet{Y} \geq (m+1)\dim Y$ and we finish by applying Theorem
\ref{thm.KLTChar}.

(3) is immediate since in the formula for $\lct(X,Y)$ the only
feature of $X$ that appears is its dimension.

The formation of jet schemes preserves products\footnote{By
definition the functor $\jet{\usc}$ has a left adjoint thus commutes
with direct products.} and hence it follows that $\dim \jet{Y \times
Y'} = \dim \jet{Y} + \dim \jet{Y'}$ from which (4) is implied
immediately.
\end{proof}

\subsection{Proof of threshold formula}
Now we present the proof of Theorem \ref{thm.KLTLC}. For this we
first recall the transformation rule in a slightly more general form
than stated above. Let $A \subseteq \jetinf{X}$ be a measurable
subset (for example a stable subset) and let $g: \jetinf{X} \to \QQ$
be a function whose level sets are measurable. Then
\[
    \int_A \LL^g d\mu_X = \int_{f_\infty^{-1}(A)} \LL^{g\circ f_{\infty}
    -\ord_{K_{X'/X}}} d\mu_{X'}
\]
provided that one of the integrals exists, which then implies the
existence of the other. The new features are minor. Clearly, it is
allowed to integrate only over a measurable subset $A$ as long as we
also only integrate over its measurable image $f^{\infty}(A)$ as
well (that this image is measurable can be deduced from the proof of
the transformation rule). In order to make the expression $\LL^{q}$
for rational $q$ defined we have to adjoin roots of $\LL$ to the
already huge ring $\Mm_k$ and define the dimension of $\LL^{q}$ as
$q$. For the following application it is in fact enough to adjoin a
single root of $\LL$ such that the new value ring of the integral is
$\Mmhat_k[\LL^{1/n}]$ for a sufficiently big $n$.

We apply this result with $g=q\cdot \ord_Y=\ord_{qY}$ and $A=\ord^{-1}_Y(m+1)$
so that $f_\infty^{-1}(A)=\ord^{-1}_F(m+1)$ (up to measure zero) and $g\circ
f_{\infty}
    -\ord_{K_{X'/X}}=-\ord_{K_{X'/X}-qF}$. Thus we get:
\begin{equation}\label{eqn.translog}
    \int_{\ord^{-1}_Y(m+1)} \LL^{q\ord_Y} d\mu_X = \int_{\ord^{-1}_F(m+1)}
    \LL^{-\ord_{K_{X'/X}-qF}} d\mu_{X'}
\end{equation}
Now, the left hand side of this equation contains information about
the dimension of the $m$th jet scheme $J_m(Y)$, whereas the right
hand side allows us to express this dimension in terms of the data
of the log resolution. Together this will lead to a proof of Theorem
\ref{thm.KLTLC}.

Recall from Section \ref{sec.measFunc} that
\[
    \ord^{-1}_Y(m+1)=\pi^{-1}_m(\jet{Y})-\pi_{m+1}^{-1}(\jet[m+1]{Y})
\]
and thus for the measure we have
\[
    \mu_X(\ord^{-1}_Y(m+1))=[\jet{Y}]\LL^{-nm}-[\jet[m+1]{Y}]\LL^{-n(m+1)}.
\]
Now computing the left hand side of equation \eqref{eqn.translog}
one gets
\[
\begin{split}
    S_m &\defeq \int_{\ord_Y(m+1)} \LL^{q\ord_Y} d\mu_X =
    \mu_X(\ord^{-1}_Y(m+1))\LL^{q(m+1)} \\
    &=([\jet{Y}]-[\jet[m+1]{Y}]\LL^{-n})\LL^{-nm+q(m+1)}
\end{split}
\]
Recalling that
\begin{equation}\label{eqn.dimcrit}
    \dim[\jet{Y}]+n \geq \dim[\jet[m+1]{Y}]
\end{equation}
this implies, while evaluating at the dimension (see section
\ref{sec.motMeasure}), that
\begin{equation}\label{eqn.A}\tag{A}
    \dim S_m \leq \dim \jet{Y} - nm + q(m+1)
\end{equation}
with ``$<$'' holding only if we have equality in
\eqref{eqn.dimcrit}. This is a consequence of the property of the
dimension function which says that $\dim (A+B) \leq \max\{\dim
A,\dim B\}$ with equality as soon as $\dim A \neq \dim B$.

Now we turn to computing the right hand side of equation
\eqref{eqn.translog}:
\[
\begin{split}
    \int_{\ord_F^{-1}(m+1)} \LL^{-\ord_{K-qF}} d\mu_{X'} &= \sum_{i=1}^\infty \mu_{X'}(
    \ord^{-1}_F(m+1) \cap \ord^{-1}_{K-qF}(i))\LL^{-i} \\
    &= \sum_{r \in A_m} \mu_{X'}(\bigcap_i \ord^{-1}_{D_i}(r_i))\LL^{-\sum
    (b_i-qa_i)r_i}
\end{split}
\]
where the last equality relies on a partitioning of $\ord_F^{-1}(m+1)$,
according to intersection with each component $D_i$ of the occurring normal
crossing divisors:
\begin{gather*}
    A_m = \{\, r=(r_1,\ldots,r_s)\,|\, r_i \geq 0, \sum a_ir_i = m+1\,\} \qquad
    \text{(this ensures $\ord_F = m+1$)} \\
    \ord^{-1}_F(m+1) = \bigsqcup_{r\in A_m} (\bigcap_i
    \ord^{-1}_{D_i}(r_i))
\end{gather*}
Clearly, this refines the partition $\ord^{-1}_F(m+1) = \bigsqcup_i
(\ord^{-1}_F(m+1) \cap \ord^{-1}_{K-qF}(i))$ and we have for $\gamma
\in \cap \ord^{-1}_{D_i}(r_i)$ that
$\ord_{K-qF}(\gamma)=\sum(b_i-qa_i)r_i$ which justifies the above
computation. The point now is that $\mu_{X'}(\bigcap_i
\ord^{-1}_{D_i}(r_i))$ was computed explicitly in Lemma
\ref{lem.compInt} to be equal to $[D^0_{\op{supp}
        r}](\LL-1)^{|\op{supp} r|}\LL^{-\sum r_i}$. Hence
\[
\begin{split}
    S_m &= \int_{\ord^{-1}_F(m+1)} \LL^{-\ord_{K-qF}} d\mu_{X'} \\
        &= \sum_{r \in A_m} [D^0_{\op{supp}
        r}](\LL-1)^{|\op{supp} r|}\LL^{-\sum(b_i-qa_i+1)r_i}
\end{split}
\]
Note that the dimension of each summand is equal to
$n-|\op{supp}r|+|\op{supp}r|-\sum(b_i-qa_i+1)r_i = n
-\sum(b_i-qa_i+1)r_i$, and since the coefficient of the highest
dimensional part of each summand has a positive sign, there is no
cancellation of highest dimensional parts in the sum. Thus we get
\begin{equation}\label{eqn.B}\tag{B}
    \dim S_m = \max_{r \in A_m}\{\, n -\sum(b_i-qa_i+1)r_i \,\}
\end{equation}
With the formulas \eqref{eqn.A} and \eqref{eqn.B} at hand the proof
of Theorem \ref{thm.KLTLC} follows easily.
\begin{proof}[Proof of Theorem \ref{thm.KLTLC}]
    It is enough to prove the first equivalence of Theorem \ref{thm.KLTLC}, the
    second one being a limiting case of the first. Slightly reformulating and using
    the definition of KLT we have to show that
\[
b_i-qa_i+1 > 0\quad \forall i \, \iff \, \dim \jet{Y}-nm+q(m+1) < n\quad
\forall m
\]
Let us first treat the implication ``$\Leftarrow$'': In fact, we
only need to assume the right hand side for one $m_0$ such that
$m_0+1$ is divisible by each $a_i$. Then using equation
\eqref{eqn.A} we have $\dim S_{m_0} \leq \dim
\jet[m_0]{Y}-nm_0+q(m_0+1) < n$. Since $b_i-qa_i+1>0$ holds
trivially for $a_i=0$ we only need to consider $i$ such that $a_i
\neq 0$. In this case (for fixed $i$) define the tuple
$r=(r_1,\ldots,r_s)$ by setting $r_i=\frac{m_0+1}{a_i}$ and $r_j=0$
otherwise. Clearly $r \in A_{m_0}$. Now, equation \eqref{eqn.B}
implies
\[
    n > \dim S_{m_0} \geq n - \sum_{j=1}^s(b_j-qa_j+1)r_j = n-(b_i-qa_i+ 1)r_i
\]
which says nothing but that $ b_i-qa_i +1>0$ as required.

Now we treat the converse ``$\Rightarrow$'': Assuming that the pair
$(X,q\cdot Y)$ is KLT (\ie $\dim S_m < n$ for all $m$ by equation
\eqref{eqn.B}) and that
\begin{equation}\label{eqn.jetineq}
    \dim \jet{Y}-nm+q(m+1) \geq n
\end{equation}
for some $m$ we seek a contradiction. These two assumptions together
imply that the inequality \eqref{eqn.A} is strict, that is we have
$\dim S_m < \dim \jet{Y}-nm +q((m+1)$. This, as we argued before,
happens only if
\[
    \dim \jet{Y} = \dim \jet[m+1]{Y} - n.
\]
Substituting this last equality into \eqref{eqn.jetineq} we get
$\eqref{eqn.jetineq}$ for $m$ replaced by $m+1$. Repeating this we
obtain
\[
    \dim \jet[m+i]{Y} = \dim \jet{Y} + in
\]
for all $i$ which contradicts Proposition \ref{prop.dimJet} below.
\end{proof}
We include here a more elementary and more explicit version of
\ref{prop.JetYzero}. This is also due to \Mustata, we only sketch
the proof here and refer the reader to \cite[Lemma
3.7]{Mustata.JetCI}.
\begin{proposition}\label{prop.dimJet}
    Let $X$ be smooth of dimension $n$ and let $Y \subseteq X$ be a subvariety.
    Let $a=\op{mult}_y Y$ be the local multiplicity of $Y$ at the point $y \in Y$, then
    \[
        \dim (\pi_Y^m)^{-1}(y) \leq m \cdot \dim X - \floor{\textstyle\frac{m}{a}}.
    \]
    Thus if $a$ is now the maximum of all local multiplicities of $Y$ one has
    \[
        \dim \jet{Y} \leq \dim Y + m \cdot \dim X - \floor{\textstyle\frac{m}{a}} \leq (m+1) \cdot \dim X - \floor{\textstyle\frac{m}{a}}.
    \]
    It follows that if $Y$ is nowhere dense in $X$, then
    $\mu_X(\jetinf{Y})=0$.
\end{proposition}
\begin{proof}
The second statement clearly follows from the first which in turn
immediately reduces to the case that $Y \subseteq X$ is a
hypersurface. Since $X$ is smooth it is \'etale over $\AA^n$ with
$y$ mapping to $0$. Since $(\pi_X^\infty)^{-1}(y)$ gets thereby
identified with $(\pi_{\AA^n}^{\infty})^{-1}(0)$ we may assume that
$Y \subseteq \AA^n$ is given by the vanishing of $f \in k[x] =
k[x_1,\ldots,x_n]$ and the point $y$ is the origin.

The condition that the local multiplicity of $Y$ at $0$ is equal to
$a$ means that the smallest degree monomial of $f$ has degree $a$ in
$x_0,\ldots,x_n$. For simplicity we assume now that $f$ is
homogeneous of degree $a$, in general one can combine the following
proof with a deformation argument to reduce to this case along the
way, see \cite{Mustata.JetCI} for this general case.

Exercise \ref{exe.EquationsJetAn} states that $\jet[m]{Y} \subseteq
\jet[m]{\AA^n}$ is given by $m+1$ equations $f^{(0)},\ldots,f^{(m)}$
in the coordinates of $\jet[m]{\AA^n}$ described in Example
\ref{ex.JetAn}. Concretely, $f^{(i)}\in k[x^{(0)},\ldots, x^{(i)}]$
is given as the coefficient of $t^i$ in the power series
\[
    f(\sum_i x_1^{(i)}t^i,\ldots,\sum_i x_n^{(i)}t^i).
\]
With the notation $f_0^{(i)} \defeq
f^{(i)}(0,x^{(1)},x^{(2)},\ldots,x^{(i)})$ (we abbreviated the
tuples $x_1^{(i)},\ldots,x_n^{(i)}$ by $x^{(i)}$) the fiber
$(\pi^m_Y)^{-1}(0)$ is given by the vanishing of
$f_0^{(1)},\ldots,f_0^{(m)}$ in the fiber $(\pi^m_X)^{-1}(0) \cong
\Spec k[x^{(1)},\ldots,x^{(m)}]$. Recall that we need to show that
the dimension of $(\pi^m_Y)^{-1}(0)$ is at most
$mn-\floor{\frac{m}{a}}$. Thus it is enough to show that the
dimension of the variety given by the vanishing of the ideal
$I_p=(f_0^{(a)},f_0^{(2a)},\ldots,f_0^{(pa)})$ is at most $(pa)n-p$.

To show this we turn to an initial ideal of $I_p$ with respect to
the degree reverse lexicographic order on
$k[x^{(1)},\ldots,x^{(pa)}]$ where the underlying ordering of the
variables $x^{(j)}_i$ is first by upper index and then by
lower.\footnote{That is we order the variables according to
$x^{(j)}_i \leq x^{(j')}_{i'}$ iff $j > j'$ or $j=j'$ and $i < i'$.
The degree reverse lexicographic order on a polynomial ring
$k[x_1,\ldots,x_n]$ with $x_1>x_2>\ldots>x_n$ is given by
\[
    A=x_1^{a_1}\cdots x_n^{a_n}>x_1^{b_1}\cdots x_n^{b_n}=B
\]
if $\deg A > \deg B$ or $\deg A=\deg B$ and $a_i<b_i$ for the last
index $k$ for which $a_i \neq b_i$. Roughly speaking, a monomial in
degree rev lex is big if it contains fewest of the cheap variables.
For example $x_2^4>x_1^3x_3>x_1^2x_2x_3$. Consult \cite[Chapter
15]{Eisenbud.CommAlg} for details.} Now it is a matter of
unravelling the definitions\footnote{Show that for two monomials $A
< B$ in $k[x]$ one has $\op{in}_<(A_0^{(j)}) <
\op{in}_<(B_0^{(j)})$. This reduces to the case that
$f=\op{in}_<(f)$ is a monomial. Then observe that the monomial
$f(x^{(j)})$ appears in $f^{(ja)}$. Since the rev-lex-ordering is
such that a term is large if it contains fewest variables which are
small it follows that $f(x^{(j)})$ is in fact the leading term as
claimed.} (of the $f^{(i)}$, of the order \dots) to see that
\[
    \op{in}_< (f^{(aj)}_0) = \op{in}_<(f)(x^{(j)}) \text{ for } j=1,\ldots, p\,.
\]
As these are $p$ many nontrivial equations in disjoint variables it
follows that the dimension of the vanishing locus of the initial
ideal of $I_p$ is at most $(pa)n-p$. Thus the same upper bound holds
for the dimension of the vanishing locus of $I_p$ itself, and a
forteriori for the dimension of $(\pi^{pa}_Y)^{-1}(0)$.
\end{proof}

\subsection{Bounds for the log canonical threshold}
In this section we show how the just derived description of the log
canonical threshold in terms of the dimension of the jet spaces lead
to some interesting bounds for $\lct(X,Y)$.
\begin{proposition}
    Let $a$ be the maximal local multiplicity of a point in $Y$. Then
    \[
        \frac{1}{a} \leq \lct(X,Y) \leq \frac{\dim X}{a}.
    \]
\end{proposition}
\begin{proof}
    We may assume that $Y$ is neither $\emptyset$ nor all of $X$
    since these cases are trivial.
    Let $p$ be a point with maximal multiplicity $a$. The second part of Exercise \ref{exe.Fiber0JnHomog}
    shows that $\dim \jet[a-1]{Y}\geq \dim (\pi^{a-1})^{-1}(p) \geq \dim X \cdot (a-1)$. Hence by
    Theorem \ref{thm.lcformula}
    \[
        \lct(X,Y) \leq \dim X - \frac{\dim \jet[a-1]{Y}}{a} \leq \dim X
        - \dim X \frac{a-1}{a} = \frac{\dim X}{a}
    \]

    For the lower bound Proposition \ref{prop.dimJet} gives
    $\dim \jet{Y} \leq \dim Y + m \cdot \dim X -
    \floor{\textstyle\frac{m}{a}}$
    and again using Theorem \ref{thm.lcformula} this yields
    \[
    \lct(X,Y)
    \geq \frac{\dim X - \dim Y}{m+1} +
    \frac{\floor{\textstyle{\frac{m}{a}}}}{m+1}
    \]
    for all $m$. For
    sufficiently divisible and sufficiently big $m$ this implies
    that $\lct(X,Y) \geq \frac{1}{a}$ as claimed.
\end{proof}

The next application of the arc space techniques is a bound for the
log canonical threshold of a homogeneous hypersurface. This uses the
following exercise as a key ingredient.
\begin{exercise}\label{exe.Fiber0JnHomog}\solution{This exercise is
solved by explicitly writing down the equations which define
$(\pi^{m}_Y)^{-1}(0)$ within $(\pi^m_{\AA^n})^{-1}(0)$. These turn
out to be the same as the defining equations of the right hand side;
just the variables are different.}
 Let $Y \subseteq \AA^n$ be a homogeneous hypersurface of degree
$d$. Show that one has an isomorphism
\[
    (\pi^{m}_Y)^{-1}(0) \cong \jet[m-d]{Y}\times \AA^{n(d-1)}
\]
for all $m \geq d-1$, where we set $\jet[-1]{Y}$ to be a point.

    Drop the assumption \emph{homogeneous} and assume instead that the
    local multiplicity of $Y$ at $p$ be equal to $a$. Show that
    $(\pi^{a-1}_Y)^{-1}(p) \cong \AA^{n(a-1)}$.
\end{exercise}
\begin{proposition}
    Let $Y \subseteq \AA^n$ be a homogeneous hypersurface of degree
    $d$. Then
    \[
        \lct(\AA^n,Y) \geq \min\left\{ \frac{n-r}{d}, 1 \right\}
    \]
    where $r = \dim \op{Sing} Y$.
\end{proposition}
\begin{proof}
    One key ingredient is the observation of the previous Exercise
    \ref{exe.Fiber0JnHomog} that for $m \geq d-1$
    \[
        (\pi^{m}_Y)^{-1}(0) \cong \jet[m-d]{Y}\times \AA^{n(d-1)}
    \]
    By semicontinuity of $\dim (\pi^m_Y)^{-1}(p)$ the inequality
    \[
        \dim (\pi^m_Y)^{-1}(p) \leq \dim \jet[m-d]{Y} + n(d-1)
    \]
    holds for all $p$, and in particular for the $p \in \op{Sing} Y$.
    Hence all together we get the estimate
    \[
        \dim \jet{Y} \leq \max\left\{ \dim \jet[m-d]{Y} + n(d-1),
        (n-1)(m+1) \right\}
    \]
    where $(n-1)(m+1)$ is equal to $\dim \jet{Y-\op{Sing} Y}$ which is
    always a lower bound for $\dim \jet{Y}$. Now set $m=pd-1$ and
    apply the inequality repeatedly to get
    \[
        \dim \jet[pd-1]{Y} \leq \max\left\{ (nd-n+r)\cdot p, (n-1)pd \right\}
    \]
    which amounts to
    \[
        \lct(\AA^n,Y) = n - \frac{\dim \jet[pd-1]{Y}}{pd} \geq \min \left\{
        \frac{n-r}{d}, 1\right\}
    \]
    since the log canonical threshold is computed via the dimension
    of the jet spaces $\jet[pd-1]{Y}$ for sufficiently divisible
    $pd$.
\end{proof}
In \cite{EinMus.LocCanThresh} the main achievement is to
characterize the extremal case as follows: In the case that
$\lct(\AA^n,Y) \neq 1$ one has $\lct(\AA^n, Y) = \frac{n-r}{d}$ if
and only if $Y \cong Y' \times \AA^r$ for some $Y' \subseteq
\AA^{n-r}$ a hypersurface.

\subsection{Inversion of adjunction}
One of the most celebrated applications of motivic integration to
birational geometry is a much improved understanding of the
Inversion of Adjunction conjecture of Shokurov
\cite{Shokurov.GeomLogFlips} and Kollar \cite{Kollar.SingPairs}. The
conjecture describes how certain invariants of singularities of
pairs behave under restriction. The following proposition goes in
this direction as it shows that under restriction to a smooth
hypersurface the log canonical threshold can only decrease, that is
the singularities cannot get better under restriction.
\begin{proposition}
    Let $(X,Y)$ be a pair and $H$ a smooth hypersurface in $X$. Then
    \[
        c_H(X,Y) \geq \lct(H,Y \cap H)
    \]
    where $c_H(X,Y)$ is the log canonical threshold of the pair $(X,Y)$ \emph{around} $H$,
    that is the minimum $\lct(U,Y \cap U)$ over all open $U \subseteq X$ with $U
    \cap H \neq \emptyset$.
\end{proposition}
\begin{proof}
    The proof uses a straightforward extension of the formula for the log canonical threshold to
    include this more general case of $c_H(X,Y)$. In fact the same
    proof as above shows that one has
    \[
        c_H(X,Y) = \dim X - \sup_m \left\{\frac{\dim_H
        \jet{Y}}{m+1}\right\}
    \]
    with equality for sufficiently divisible $m+1$. Here
    $\dim_H \jet{Y}$ is the \emph{dimension of $\jet{Y}$
    along $H$}, that is the maximal dimension of an irreducible
    component $T$ of $\jet{Y}$ such that $\pi^m(T) \cap H \neq \emptyset$.

    Let $T$ be an irreducible component of $\jet{Y}$ such that $\pi^m(T) \cap H \neq \emptyset$.
    This implies that $T \cap \jet{Y \cap H}$ is also nonempty since
    the projection $\pi^m_{Y\cap H}$ is surjective. Since $H \subseteq X$
    is locally given by one equation the same is true for $H \cap Y
    \subseteq Y$. Hence $\jet{H \cap Y} \subseteq \jet{Y}$ is given
    by at most $m+1$ equations (see Exercise
    \ref{exe.EquationsJetAn}). Hence $\dim_H \jet{Y} \leq \dim
    \jet{Y\cap H} +(m+1)$ which shows
    \[
        \sup_m \left\{\frac{\dim_H \jet{Y}}{m+1}\right\} \leq \sup_m \left\{\frac{\dim_H \jet{Y \cap
        H}}{m+1}\right\}+1
    \]
    which implies the claimed inequality $c_H(X,Y) \geq
    \lct(H,Y \cap H)$.
\end{proof}

The \emph{Inversion of Adjunction} Conjecture of Koll\'ar and
Shokurov describes how the singularities of pair behave under
restriction to a Cartier divisor. More precisely, let $(X,Y)$ be a
pair where we allow $Y=\sum q_i Y_i$ to be any formal integer
combination (rational or real combination even) of closed
subvarieties of $X$. With the notation $f:X' \to X$ of a log
resolution as above (in particular $f^{-1}Y=\sum a_i E_i$ and
$K_{X'/X}=\sum b_iE_i$) and a subvariety $W \subseteq X$ fixed we
define the \emph{minimal log discrepancy}
\[
    \mld(W;X,Y) \defeq \begin{cases}
                              \min\{b_i-a_i+1 | f(E_i) \subseteq W \}& \text{if this minimum is
                              non-negative}\\
                              - \infty &\text{otherwise.}
                        \end{cases}
\]
It follows that $(X,Y)$ is log canonical on an open subset
containing $W$ iff $\mld(W;X,Y) \neq -\infty$. The inversion of
adjunction conjecture now states:
\begin{conjecture}
    With $(X,Y)$ and as above, let $D$ be a normal effective Cartier
    divisor on $X$ such that $D \not\subseteq Y$ and let $W \subset D$ a proper closed
    subset. Then we have
    \[
        \mld(W;X,Y+D) = \mld(W;D,Y|_D).
    \]
\end{conjecture}
The inequality ``$\leq$'' is the \emph{adjunction} part and is well
known to follow from the adjunction formula $K_D = (K_X +D)|_D$. The
reverse inequality ``$\geq$'' is the critical part of this
conjecture.

In \cite{EinMustYas.JetDiscr} the conjecture was proved in the case
that $X$ is smooth and $Y$ is effective. In \cite{EinMus.InvAdj}
this was established even for $X$ a complete intersection (and $Y$
effective). The proof of these results use a description of the
minimal log discrepancies in terms of dimensions of certain
cylinders of the jet spaces of $X$, analogous to the one for the log
canonical threshold.
\begin{proposition}
With the notation as above and for $X$ smooth and $Y$ effective,
\begin{gather*}
    \mld(W;X,Y)  \geq \tau \\  \iff   \\
    \qquad \codim_{\jetinf{X}}( \Cont^\nu_Y \cap \pi_0^{-1}(W)) \geq  \sum q_i\nu_i + \tau \qquad \text{ for all
    multi-indices $\nu$.}
\end{gather*}
where $\Cont^\nu_Y = \cap_i \ord^{-1}_{Y_i}(\nu_i)$.
\end{proposition}

The proof of this is not more than a technical complication of the
proof of the log canonical threshold formula we gave above. With
this characterization of $\mld(W;X,Y)$ the proof of inversion of
adjunction becomes a matter of determining the co-dimensions of the
cylinders involve. In the case that $X$ and $D$ are both smooth this
is quite easy (and could be done as an exercise). In general ($X$ a
complete intersection) the combinatorics involved can become quite
intricate, \cf \cite{EinMus.InvAdj}.

One should point out that after these results on inversion of
adjunction were obtained, Takagi \cite{Takagi.InvAdj} found an
alternative approach using positive characteristic methods.

\subsection{Geometry of arc spaces without explicit motivic
integration.} As it should have become apparent by now, the
applications of motivic integration to birational geometry are by
means of describing certain properties of a variety $X$ in terms of
(mostly simpler) properties of its jet spaces $\jet{X}$. Motivic
integration serves as the path to make this connection. However, due
to some combinatiorial difficulties one encounters along this path,
one can ask if there is a more direct relationship. This is indeed
the case and it is the content of the paper of
\cite{EinLazMus.ContLoci} of Ein, Lazarsfeld and \Mustata which is
the source of the material in this section.

Their point is that instead of using the birational transformation
rule to control the dimension of components of the jet spaces one
uses the Key Theorem \ref{thm.maintechnical} of its proof to more
directly get to the desired information.

We keep the notation of a subvariety $Y \subseteq X$ of a smooth
variety $X$ and define.
\begin{definition}
Let $Y$ be a subvariety of $X$ and $p \geq 0$ an integer define the
\emph{contact locus} to be the cylinder
\[
    \Cont^p_Y \defeq \ord^{-1}_Y(p) \subseteq  \jetinf{X}
\]
\end{definition}
The aim is to understand the components (or at least the dimension)
of the cylinders $\Cont^p_Y$ in terms of the data coming from a log
resolution of the pair $(X,Y)$. Using the notation of Section
\ref{sec.notbir} we fix a log resolution $f:X' \to X$ of the pair
$(X,Y)$ and denote $f^{-1}Y = \sum_1^k a_iE_i$ and $K_{X'/X} =
\sum_1^k b_i E_i$ where the support of $\sum E_i$ is a simple normal
crossing divisor.

\begin{definition}
    Given $E = \sum_1^k E_i$, a simple normal crossing divisor of $X'$, and
    a multi-index $\nu = (\nu_1,\ldots,\nu_k)$ define the \emph{multi
    contact locus}
    \[
        \Cont^\nu_E \defeq \left\{\, \gamma' \in \jetinf{X'}\, |\,
        \ord_{E_i}(\gamma') = \nu_i \text{ for }i=1\ldots
        k\,\right\}.
    \]
\end{definition}
\begin{definition}
    For every cylinder $C \subseteq \jetinf{X}$ there is a well defined
    notion of codimension, namely
    \[
        \codim C \defeq \codim(\jet{X},\pi_m C)
    \]
    for $m \gg 0$.
\end{definition}
Of course one must check that this is independent of the chosen $m
\gg 0$. This however is immediately clear from the definition of
cylinder.
\begin{exercise}\label{exe.CodimContNCD}
    For $C \subseteq \jetinf{X}$ a cylinder, show that $\codim C =
    \dim X - \dim \mu_X (C)$.
\end{exercise}
The following proposition replaces in this new setup the computation
of the motivic integral of a normal crossing divisor in Proposition
\ref{prop.comp}. Note the comparative simplicity!
\begin{proposition}
    For $E = \sum_1^k E_i$ a simple normal crossing divisor $\Cont^\nu_E$ is a
    smooth irreducible cylinder of codimension
    \[
        \codim \Cont^\nu_E = \sum_1^k \nu_i
    \]
    provided $\Cont^\nu_E$ is nonempty.
\end{proposition}
\begin{proof}
    This is a computation in local coordinates. Assume $E$ is
    locally given by the vanishing of the first $k$ of the
    coordinates $x_1=\ldots=x_k=0$. Then, for an arc $\gamma$, which
    is determined by $\gamma(x_i)=\sum \gamma_i^{(j)}t^j$ for $i=1\ldots n$, to have the
    prescribed contact order with the $E_i$ means precisely that
    $\gamma_i^{(j)} = 0$ for $j< \nu_i$, and $\gamma_i^{(\nu_i)} \neq
    0$. Hence for $m \gg 0$ we have
    \[
        \pi_m(\Cont^\nu_E) \cong (\AA^1-\{0\})^n \times
    \prod\AA^{m-\nu_i}
    \]
    and therefore $\Cont^\nu_E$ is smooth
    irreducible and of codimension $\sum \nu_i$.
\end{proof}
The central result (replacing the transformation rule) is the
following Theorem
\begin{theorem}\label{thm.ContactLociImages}
    With the notation as above one has for all $p > 0$ a finite partition
    \[
        \Cont^p_Y = \bigsqcup_v f_\infty \Cont^n_E
    \]
    where the disjoint union is over all multi-indices $\nu$ such
    that $\sum \nu_i a_i = p$.

    For every multi-index is the set $f_\infty \Cont^v_E$ an
    irreducible cylinder of codimension \[\sum \nu_i(b_i+1).\]

    In particular, for each irreducible component $Z$ of\/ $\Cont^p_Y$
    there is a unique multi-index $\nu$ such that $\Cont^\nu_E$
    dominates $Z$.
\end{theorem}
\begin{proof}
    As in the proof of the transformation rule the key ingredient is
    Theorem \ref{thm.maintechnical}. With this at hand the proof is
    not difficult.

    The condition $\sum \nu_i a_i = p$ ensures that $f_\infty\Cont^\nu_E \subseteq
    \Cont^p_Y$. The surjectivity of $f_\infty$ (\cf Exercise \ref{exe.finftysurj}) on the other hand
    implies the reverse inclusion. The disjoined-ness of the union follows from the
    fact that there is a one-to-one map between
    \[
        \jetinf{X'}-\jetinf{E} \to[1-1] \jetinf{X}-\jetinf{Y}
    \]
    (which is an implication of the valuative criterion for properness as explained in Section \ref{sec.Finfty})
    and the observation that each $\Cont^\nu_E$ is contained in the left hand
    side.

    Since $\Cont^\nu_E \subseteq \Cont^{\sum b_iv_i}_{K_{X'/X}}=\ord^{-1}_{K_{X'/X}}$ it
    follows from Theorem \ref{thm.maintechnical} (a) and Proposition \ref{prop.imageCyl} that $f_\infty\Cont^\nu_E$
    is a cylinder. Part (b) of Theorem \ref{thm.maintechnical} shows
    that
    \[
        f_m : \Cont^\nu_E \to f_m\Cont^\nu_E
    \]
    is a piecewise trivial $\AA^{\sum b_i\nu_i}$--fibration. By
    Exercise \ref{exe.CodimContNCD} is the codimension of
    $\Cont^\nu_E$ equal to $\sum \nu_i$. Hence the codimension of
    its image under $f_\infty$ is
    \[
        \codim (f_\infty \Cont^\nu_E) = \sum\nu_i + \sum \nu_ib_i =
        \sum \nu_i(b_i+1)
    \]
    as claimed.
\end{proof}
As an illustration of this result we recover a very clean proof of
the log canonical threshold formula of Theorem \ref{thm.lcformula}.

\begin{proof}[Proof of Theorem \ref{thm.lcformula}]
    Let $V_m$ be an irreducible component of $\jet[m]{Y}$. For some
    $p \geq m+1$ the set $\Cont^p_Y \cap (\pi^X_m)^{-1}V_m$ is open
    in $(\pi^X_m)^{-1}V_m$, namely for the smallest $p$ (automatically $\geq m+1$ since
    each arc in $V_m$ has contact order $\geq m+1$ with $Y$) such that $\Cont^p_Y \cap
    (\pi^X_m)^{-1}V_m \neq \emptyset$. Hence there is an irreducible
    component $W$ of $\Cont^p_Y$ such that the closure
    $\overline{W}$ contains $(\pi^X_m)^{-1}V_m$.
    By Theorem \ref{thm.ContactLociImages} there is a
    unique multi-index $\nu$ (necessarily $\sum \nu_ib_i = p$) such that $\Cont^\nu_E$ dominates $W$.

    By definition of the log canonical threshold we have $b_i+1 \geq
    \lct(X,Y)a_i$ for all $i$ such that we obtain the following
    inequalities:
    \begin{equation*}
    \begin{split}
        \codim (V_m,\jet[m]{X})
                &\geq \codim W \\
                &= \codim f_\infty \Cont^\nu_E \\
                &\geq \sum \nu_i(b_i+1) \\
                &\geq \sum \nu_i \lct(X,Y)a_i \\
                &= \lct(X,Y) \cdot p = \lct(X,Y) \cdot (m+1)
    \end{split}
    \end{equation*}
    As this holds for every irreducible component of $\jet[m]{Y}$ we
    get
    \[
        \lct(X,Y) \leq \frac{\codim(\jet{Y},\jet{X})}{m+1}.
    \]
    To see that there is equality for some $m$ we pick an index $i$ such that
    $\lct(X,Y) = \frac{b_i+1}{a_i}$ and $m+1$ divisible by $a_i$. Let $\nu$ be the multi-index
    which is zero everywhere except at the $i$th spot, where it is
    $\frac{m+1}{a_i}$. Then $f_\infty \Cont^\nu_E \subseteq
    \Cont^{m+1}_Y \subseteq (\pi^X_m)^{-1} \jet{Y}$ and by Theorem
    \ref{thm.ContactLociImages} the codimension $f_\infty
    \Cont^\nu_E$ is equal to $\frac{m+1}{a_i}(b_i+1) =
    \lct(X,Y)\cdot (m+1)$. Hence in particular
    $\codim(\jet{Y},\jet{X}) \leq \lct(X,Y) \cdot (m+1)$ for this chosen $m+1$. This finishes the argument.
\end{proof}
In summary, the above agument shows that the irreducible components
$V$ of $\jet{Y}$ of maximal possible dimension, that is the ones
that compute the log canonical theshold as
$\lct(X,Y)=\codim(\jet{X},V)$ are dominated by multi-contact loci
$\Cont^\nu_E$ with $\nu_i \neq 0$ for all the indices $i$ such that
$E_i$ computes the the log canonical theshold (meaning
$\lct(X,Y)=\frac{b_i+1}{a_i}$).

We want to finish these notes with \Mustata's characterization of
rational singularities for complete intersections in terms of arc
spaces. This was indeed the first application of motivic integration
to characterizing singularities. With the just developed viewpoint
this result is not too difficult anymore.

\begin{theorem}
    Let $Y \subseteq X$ be a reduced and irreducible locally
    complete intersection subvariety of codimension $c$. Then the
    jet spaces $\jet{Y}$ are irreducible for all $m$ if and only if
    $Y$ has rational singularities.
\end{theorem}
\begin{proof}
    Let $f: X' \to X$ be a log resolution of $(X,Y)$ which dominates
    the blowup of $X$ along $Y$. Keeping the previous notation we
    may assume that $E_1$ is the exceptional divisor of this blowup.
    In \cite{Must.JetLCI} Theorem 2.1 it is shown that $Y$ has at
    worst rational singularities if and only if $b_i \geq c a_i$ for
    every $i \geq 2$. Hence we must show
    \[
        \jet{Y} \text{ is irreducible for all $m \geq 1$} \iff b_i
        \geq c a_i \text{ for $i \geq 2$}
    \]
    Assume that $\jet{Y}$ is not irreducible, that is we have a
    component $V \subseteq \jet{Y}$ other than the main component
    $\jet{Y-\Sing Y}$. As in the previous proof we have $W \subseteq
    \Cont^p_Y$ with $p \geq m+1$ whose closure contains $\pi_m^{-1}(V)$.
    By Theorem \ref{thm.ContactLociImages} this component is
    dominated by some multi-contact locus $\Cont^\nu_E$ for $\nu
    \neq (m+1,0,\ldots,0)$ since the latter is the multi-index
    corresponding to the multi-contact locus dominating $\pi_m^{-1}(\jet{Y-\Sing
    Y})$. Since $Y \subseteq X$ is a local complete intersection of
    codimension $c$ we have $\codim(V,\jet{X}) \leq (m+1)\cdot c$.
    To arrive at a contradiction assume now that $Y$ has rational
    singularities, that is assume that $b_i
        \geq c a_i$ for $i \geq 2$. Then
    \[
    \begin{split}
        (m+1) \cdot c &\geq  \codim(W) \\
                      &= \nu_1\cdot c \sum_{i\geq 2} \nu_i(b_i+1) \\
                      &\geq c \cdot \sum_{i\geq 1} \nu_i a_i +
                      \sum_{\geq 2} \nu_i \qquad \text{ (since $b_i
        \geq c a_i$)}\\
                      &= c \cdot p + \sum_{i \geq 2} \nu_i \geq c
                      \cdot (m+1) + \sum_{i \geq 2} \nu_i
    \end{split}
    \]
    Hence for $i \geq 2$ we must have $\nu_i = 0$, a contradiction.

    Conversely, suppose $b_i < c \cdot a_i$ for some $i \geq 2$.
    Setting $v$ to be the multi-index with all entries zero except
    the $i$th equal to $1$. Let $(m+1)=a_i$, then the $\Cont^\nu_E$
    maps to an irreducible subset $W \subseteq \Cont^{m+1}_Y$ of
    codimension $(m+1)\cdot c$. Hence $\pi_{m+1}(W)$ is an
    irreducible component of $\jet{Y}$ of codimension $(m+1)\cdot c$
    which is not the component $\overline{\jet{Y-\Sing Y}}$. Hence
    $\jet{Y}$ is not irreducible.
\end{proof}

\appendix

\section{An elementary proof of the Transformation rule.}
We present Looijenga's \cite{Looijenga.motivicmeasure} elementary
proof of Theorem \ref{thm.maintechnical} which then leads to a proof
of the transformation formula avoiding weak factorization. For this
we have to investigate more carefully the definition of the relative
canonical divisor $K_{X'/X}$ and suitably interpret the contact
multiplicity of an arc $\gamma$ with $K_{X/X'}$.

\subsection{The relative canonical divisor and differentials}
Let us consider the first fundamental exact sequence for K\"ahler
differentials, as it plays a pivotal role in all that follows.

The morphism $f: X' \to X$ induces a linear map, its derivative,
$f^*\Omega_X \to[df] \Omega_{X'}$ which is part of the first
fundamental exact sequence for K\"ahler differentials:
\begin{equation}\label{eq.fundexdeR}
     0 \to f^*\Omega_X \to[df] \Omega_{X'} \to \Omega_{X'/X} \to 0
\end{equation}
Note that by our assumption of smoothness, the $\Oo_{X'}$--modules
$f^*\Omega_X$ and $\Omega_{X'}$ are locally free of rank $n=\dim X$.
Since, by birationality of $f$, $\Omega_{X'/X}$ has rank zero, the
first map is injective as well. Taking the $n$th exterior power we
obtain the map
\[
    0 \to f^*\Omega_X^n \to[\wedge^n df] \Omega^n_{X'}
\]
of locally free $\Oo_{X'}$ modules of rank $1$. If we set $\omega =
\Omega^n$ and tensor the above sequence with the invertible sheaf
$\omega_{X'}^{-1}$ we obtain
\[
    f^*\omega_X \tensor \omega_{X'}^{-1} \subseteq \Oo_{X'}
\]
thus identifying $f^*\omega_X \tensor \omega_{X'}^{-1}$ with a
locally principal ideal in $\Oo_{X'}$, which we shall denote by
$J_{X'/X}$ (so, by definition, $J_{X'/X}$ is the 0-th Fitting ideal
of $\Omega_{X'/X}$). Now define $K_{X'/X}$ to be the Cartier divisor
which is locally given by the vanishing of $J_{X'/X}$. It is
important to note that $K_{X'/X}$ is defined as an effective divisor
and not just as a divisor class. By choosing bases for the free
$\Oo_{X'}$-modules $f^*\Omega_X$ and $\Omega_{X'}$ the map $df$ is
given by a $n \times n$ matrix with entries in $\Oo_{X'}$. Its
determinant is a local defining equation for $K_{X'/X}$.

Let $L$ be an extension field of $k$ and let $\gamma: \Spec L\lbr t
\rbr \to X$ be a $L$-rational point of $\jetinf{X'}$, and assume
that $\ord_{K_{X'/X}}(\gamma)=e$. By definition of contact order,
this means that $(t^e) = \gamma^*(J_{X'/X}) \subseteq L\lbr t \rbr$.
As $J_{X'/X}$ is locally generated by $\det df \in \Oo_{X'}$ (well
defined up to unit) we obtain that $(t^e)= \det(\gamma^*(df))$. The
pullback of the sequence \eqref{eq.fundexdeR} along $\gamma$
illustrates the situation:
\begin{equation}\label{eqn.deRgamma}
    0 \to (f\circ\gamma)^*\Omega_X \to[\gamma^*df] \gamma^*\Omega_{X'} \to \gamma^*\Omega_{X'/X} \to 0
\end{equation}
Since $L\lbr t \rbr$ is a PID, we can choose bases of
$(f\circ\gamma)^*\Omega_X$ and $\gamma^*\Omega_{X'}$ such that
$\gamma^*(df)$, a map of free $L\lbr t \rbr$ modules of rank $n$, is
given by a diagonal matrix. With respect to this basis the exact
sequence \eqref{eq.fundexdeR} takes the form
\begin{equation}\label{eqn.deRgammaBasis}
    0 \to L\lbr t \rbr^n \to[{\left(\begin{smallmatrix} t^{e_1}&&0 \\ &\ddots& \\ 0&&t^{e_n} \end{smallmatrix}\right)}]
    L\lbr t \rbr^n \to \oplus \frac{L \lbr t \rbr}{(t^{e_i})} \to 0
\end{equation}
The condition that $\ord_{K_{X'/X}}(\gamma)=e$ translates into
$\sum_{i=1}^n e_i=e$ or, equivalently, into saying that the
rightmost module is torsion of length $e$.

\subsection{Proof of Theorem \ref{thm.maintechnical}}\label{sec.proofTrans} We start by recalling the
statement of Theorem \ref{thm.maintechnical} we want to proof
slightly reformulated in order to set up the notation that is used
in its proof below.

\begin{theorem}
    Let $f:{X^\prime} \to X$ be a proper birational morphism of smooth
    $k$-varieties. Let $C_e'= \ord_{K_{X'/X}}^{-1}(e)$ where $K_{X'/X}$ is
    the relative canonical divisor and let $C_e \defeq f_\infty C_e'$.
    Let $\gamma \in C_e'$  an $L$-point of $\Jj_\infty({X^\prime})$,
    that is a map $\gamma^*: \Oo_{X^\prime} \to L \lbr t \rbr $, satisfying $\gamma^*(J_{X'/X})=(t^e)$,
    with $L \supseteq k$ a field extension. Then for $m \geq 2e$ one has:

    \begin{enumerate}
    \item[(a')] For all $\xi \in \Jj_\infty(X)$ such that
    $\pi^X_m(\xi)=f_m(\pi^{X'}_m(\gamma))$ there is $\gamma' \in \Jj_\infty({X^\prime})$ such
    that $f_\infty(\gamma')=\xi$ and $\pi^{X'}_{m-e}(\gamma')=\pi^{X'}_{m-e}(\gamma)$.
    In particular, the fiber of $f_m$ over $f_m(\gamma_m)$ lies in the fiber of
    $\pi^m_{m-e}$ over $\gamma_{m-e}$.
    \item[(a)] $\pi_m(C_e')$ is a union of fibers of
    $f_m$.\xfootnote{In an earlier version of the paper we had here: ``$C_e=(\pi^X_m)^{-1}(f_m(\pi^{X'}_m(C_e')))$
    and is hence stable at level $m$.'' This is also implied by (a'): Since $f_m\pi^{X'}_m=\pi^X_m f_\infty$ it is
    clear that ``$\subseteq$''
holds. To check ``$\supseteq$'' recall that $C_e'$ is stable at
level $\geq e$ and since $m\geq 2e$ one has
$C_e'=(\pi^{X'}_{m-e})^{-1}\pi^{X'}_{m-e}C_e'$. Now (a') shows that
\[
\begin{split}
    (\pi^X_m)^{-1}\pi^X_m f_\infty(\ord^{-1}_K(e)) &\subseteq
    f_{\infty}((\pi^{X'}_{m-e})^{-1}\pi^{X'}_{m-e}(\ord^{-1}_K(e))) \\
        &=f_\infty(\ord^{-1}_K(e)).
\end{split}
\]
To justify the inclusion above, let $\xi \in (\pi^X_m)^{-1}\pi^X_m
f_\infty(\ord^{-1}_K(e)) = (\pi^X_m)^{-1} f_m \pi^{X'}_m
(\ord^{-1}_K(e))$; in particular $\pi^X_m(\xi)=
f_m\pi^{X'}_m(\gamma)$ for some $\gamma \in C_e'$. Then part (a')
yields a $\gamma' \in \jetinf{X}$ such that
$\pi^{X'}_{m-e}(\gamma')=\pi^{X'}_{m-e}(\gamma')\in
\pi^{X'}_{m-e}(\ord^{-1}_K(e))$ and $f_\infty(\gamma') = \xi$ as
claimed.}
    \item[(b)] The map $f_m: \pi^{X'}_m(C_e') \to C_e$ is a piecewise
    trivial $\AA^e$ fibration.
    \end{enumerate}
\end{theorem}
\begin{proof}
To ease notation we will denote truncation by lower index, i.e.\
write $\gamma_m$ as shorthand for $\pi^X_m(\gamma)$. We already
pointed out before that (a') implies (a):

The proof of (b) can be divided into two steps.  First we show that
the fiber of $f_m$ over $f_m(\gamma_m)$ can be naturally identified
with $\Der_{\Oo_X}(\Oo_{X^\prime},\frac{L \lbr t \rbr
}{(t^{m+1})})$. Then we have to show that the latter is an affine
space of dimension $e$. As this is easy let's do it first:
Immediately preceding this proposition we noted that the cokernel of
$\gamma^*(df)$ is torsion of length $e$ as a $L \lbr t \rbr$-module.
This cokernel is $\gamma^*\Omega_{X'/X}$. Since $m>e$ the dual,
$\Hom_{L\lbr t \rbr}(\gamma^*\Omega_{X'/X},\frac{L\lbr
t\rbr}{(t^{m+1})})$, is also torsion of length $e$. Using
adjointness of $\gamma^*$ and $\gamma_*$ this $\Hom$ is just
$\Hom_{\Oo_{X'}}(\Omega_{X'/X},\gamma_*\frac{L\lbr
t\rbr}{(t^{m+1})})$, which is equal to
$\Der_{\Oo_X}(\Oo_{X^\prime},\gamma_*\frac{L \lbr t \rbr
}{(t^{m+1})})$ essentially by definition of $\Omega_{X'/X}$. This
shows that $\Der_{\Oo_X}(\Oo_{X^\prime},\gamma_*\frac{L \lbr t \rbr
}{(t^{m+1})})$ is isomorphic to $\AA^e_L$. Thus we are left to show
the identification $(\diamondsuit\diamondsuit\diamondsuit)$ of the
following diagram the last line of which is the first exact sequence
for derivations, analogous to the above exact sequence of K\"ahler
differentials.
\[
\xymatrix{
    {\left(\txt{fiber of $f_m$ \\ over $f_m(\gamma_m)$}\right)} \ar@{^(->}^{\text{(a')}}[r] \ar@{=}^{(\diamondsuit\diamondsuit\diamondsuit)}[dd] &{\left(\txt{fiber of $\pi^m_{m-e}$\\ over $\gamma_{m-e}$}\right)} \ar^{(\diamondsuit)}@{=}[d] \\
    {} &{\Der_L(\Oo_{X^\prime},\frac{(t^{m+1-e})}{(t^{m+1})})} \ar@{^(->}^{\text{induced from inclusion }\frac{t^{m+1-e}}{t^{m+1}} \subseteq \frac{L \lbr t \rbr }{t^{m+1}}}[d]  \\
    {\Der_{\Oo_X}(\Oo_{X^\prime},\frac{L \lbr t \rbr }{(t^{m+1})})} \ar@{^{(}->}[r]\ar@{^(->}^{(\diamondsuit\diamondsuit)}[ur] &{\Der_L(\Oo_{X^\prime},\frac{L \lbr t \rbr }{(t^{m+1})})} \ar[r] &\Der_L(f^*\Oo_X,\frac{L \lbr t \rbr }{(t^{m+1})}) \\
    }
\]

The identification $(\diamondsuit)$ is given by sending $\gamma'_m$
to $\gamma_m'-\gamma_m$ which, since $m \geq 2e$, can
easily\footnote{Fix an homomorphism $\gamma: R \to S$ of
$k$-algebras which makes $S$ into an $R$-algebra. For any ideal $I
\in S$ with $I^2=0$ one has a map
\begin{equation}\label{eqn.der}
    \{\,\gamma' \in \Hom_{k-\text{alg}}(R,S)\,|\,\Image(\gamma'-\gamma)\subseteq I\,\} \to \Der_k(R,I)
\end{equation}
by sending $\gamma'$ to $\gamma'-\gamma$. To check that
$(\gamma'-\gamma)$ is indeed a derivation one has to make the
following calculation verifying the Leibniz rule (note that the $R$
algebra structure on $S$ is via $\gamma$):
\[
\begin{split}
    (\gamma'-\gamma)(xy)&-((\gamma'-\gamma)(x)\gamma(y)+\gamma(x)(\gamma'-\gamma)(y))
    \\
    &=\gamma'(x)\gamma'(y)-\gamma(x)\gamma(y)-\gamma'(x)\gamma(y)+\gamma(x)\gamma(y)-\gamma(x)\gamma'(y)+\gamma(x)\gamma(y)
    \\ &=\gamma'(x)(\gamma'(y)-\gamma(y))-\gamma(x)(\gamma'(y)-\gamma(y)) \\
    &= (\gamma'-\gamma)(x)\cdot(\gamma'-\gamma)(y) = 0
\end{split}
\]
The last line is zero by the assumption that $\Image(\gamma'-\gamma)
\subseteq I$ and $I^2=0$. The obvious inverse map sending a
derivation $\delta$ to $\gamma +\delta$ shows that the two sets in
\eqref{eqn.der} are equal. This setup clearly applies in our
situation: $R=\Oo_{X^\prime}$, $S =L\lbr t \rbr / t^{m+1}$,
$I=(t^{m+1-e})$ where $m \geq 2e$ ensures that $I^2=0$.} be checked
to define an $L$-derivation $\Oo_{X^\prime} \to
\frac{(t^{m+1-e})}{(t^{m+1})}$. In this way (and using (a')) we
think of $f_m^{-1}(f_m(\gamma_m))$ as a subspace of
$\Der_L(\Oo_{X^\prime},\frac{(t^{m+1-e})}{(t^{m+1})})$. As this is
the $t^e$-torsion part of $\Der_L(\Oo_{X^\prime},\frac{L \lbr t \rbr
}{(t^{m+1})})$ and since we just observed that
$\Der_{\Oo_X}(\Oo_{X^\prime},\frac{L \lbr t \rbr }{(t^{m+1})})$ is
torsion of lenght $e$ the inclusion $(\diamondsuit\diamondsuit)$ is
also clear, and thus $(\diamondsuit\diamondsuit\diamondsuit)$
becomes a statement about subsets of
$\Der_L(\Oo_{X^\prime},\frac{(t^{m+1-e})}{(t^{m+1})})$.

Let $(\gamma'_m-\gamma_m) \in
\Der_L(\Oo_{X^\prime},\frac{(t^{m+1-e})}{(t^{m+1})})$. The image of
$(\gamma'_m-\gamma_m)$ in $\Der_L(f^*\Oo_X,\frac{L \lbr t \rbr
}{(t^{m+1})})$ is $\gamma'_m \circ f-\gamma_m \circ f$. This is zero
(i.e. $\gamma'_m \in \Der_{\Oo_X}(\Oo_{X^\prime},\frac{L \lbr t \rbr
}{(t^{m+1})})$) if and only if $f_m(\gamma'_m)= f_m(\gamma_m)$, that
is if and only if $\gamma'_m$ is in the fiber of $f_m$ over
$f_m(\gamma_m)$. This concludes the proof of (b).

In order to come by the element $\gamma' \in \jetinf{X'}$ as claimed
in (a') we construct a sequence of arcs $\gamma^k \in \jetinf{X'}$
satisfying the following two properties for all $k \geq m$:
\begin{enumerate}
    \item $\pi_{k}(f_\infty(\gamma^k)) = \pi_k(\xi)$ and
    \item $\pi_{k-1-e}(\gamma^k) = \pi_{k-1-e}(\gamma^{k-1})$ and $\pi_{m-e}(\gamma^k) = \pi_{m-e}(\gamma)$.
\end{enumerate}
Clearly, setting $\gamma^{m-1}=\gamma^m=\gamma$ these conditions
hold for $k=m$. Furthermore, the second condition implies that the
limit $\gamma' \defeq \lim_k \gamma^k$ exists and that
$\pi_{m-e}(\gamma')=\pi_{m-e}(\gamma)$. The first condition shows
that $f_\infty(\gamma')=\xi$. Thus we are left with constructing the
sequence $\gamma^k$. This is done inductively. As we already
verified the solution for $k=m$ we assume to have $\gamma^k$ and
$\gamma^{k-1}$ as claimed -- now $\gamma^{k+1}$ is constructed as
follows:

Since $\pi_k(f(\gamma^k)) = \pi_k(\xi)$ we can view their difference
as a derivation $\delta = \xi  - f \circ \gamma^k \in
\Der_{L}(\Oo_{X},\frac{(t^{k+1})}{(t^{k+2})})$ which we identify
with $\Hom_{L\lbr t
\rbr}(\gamma^{k*}f^*\Omega_X,\frac{(t^{k+1})}{(t^{k+2})})$. The
latter module appears in $\Hom_{L\lbr t\rbr}( \usc,\frac{L\lbr t
\rbr}{(t^{k+2})})$ applied to the sequence \eqref{eqn.deRgamma},
where $\gamma^k$ takes the place of $\gamma$:
\begin{equation}
\xymatrix@C=0.1pc@R=1pc{
    {\Hom(\gamma^{k*}\Omega_{X/X'},\frac{L\lbr t \rbr}{(t^{k+2})})}
    \ar@{^(->}[r] &{\Hom(\gamma^{k*}\Omega_{X'},\frac{L\lbr t \rbr}{(t^{k+2})})} \ar^{df}[rr]
    &&{\Hom(\gamma^{k*}f^*\Omega_X,\frac{L\lbr t \rbr}{(t^{k+2})})}
    \\
    {} & {\Hom(\gamma^{k*}\Omega_{X'},\frac{(t^{k+1-e})}{(t^{k+2})})\ni \delta'}\ar@{|->}[rr]\ar@{}|\bigcup[u]
    &&
    {\delta \in
    \Hom(\gamma^{k*}f^*\Omega_X,\frac{(t^{k+1})}{(t^{k+2})})}\ar@{}|\bigcup[u]}
\end{equation}
In order to understand this better we turn to the same sequence, but
with respect to the basis as in sequence \eqref{eqn.deRgammaBasis},
where it takes this form:
\[
\xymatrix@C=2pc@R=1pc{
    {\oplus \frac{(t^{k+1-e_i})}{(t^{k+2})}} \ar@{^(->}[r] &{(\frac{L\lbr t \rbr}{(t^{k+2})})^n} \ar^{\left(\begin{smallmatrix} t^{e_1}&&0 \\ &\ddots& \\ 0&&t^{e_n}
    \end{smallmatrix}\right)}[rr] && {(\frac{L\lbr t \rbr}{(t^{k+2})})^n} \\
    & {(\frac{(t^{k+1-e})}{(t^{k+2})})^n}\ar@{}|\bigcup[u] && {\delta \in
    (\frac{t^{k+1}}{(t^{k+2})})^n\ar@{}|\bigcup[u]}
}
\]

Now it becomes clear that $\delta$ lies in the image of $df$ since
$e$ and therefore all $e_i$ are less than $m+1 \leq k+1$.
Furthermore, any pre-image $\delta'$ must lie in
$\Hom(\gamma^{k*}\Omega_{X'},\frac{(t^{k+1-e})}{(t^{k+2})})$ by the
shape of the matrix and the fact that $e_i \leq e$ for all $e$. Now
pick any such pre-image $\delta'$ and define $\gamma^{k+1} \defeq
\delta' + \gamma^k$. This is an arc in $X'$ with
$\pi_{k-e}(\gamma^{k+1}) = \pi_{k-e}(\gamma^{k})$. Furthermore since
$df(\delta') = \delta$ we get
\[
     \gamma^{k+1} \circ f - \gamma^k \circ f = \delta = \xi - \gamma^k \circ f \mod
     (t^{k+2})
\]
and thus $\pi_{k+1}(\gamma^{k+1} \circ f) =
\pi_{k+1}(f_\infty(\gamma^{k+1})) = \pi_{k+1}(\xi)$.
\end{proof}
With this proof of Theorem \ref{thm.maintechnical} at hand a proof
of the Transformation rule follows immediately as indicated in
Section \ref{sec.trans.rule}.

\providecommand{\bysame}{\leavevmode\hbox
to3em{\hrulefill}\thinspace}

\end{document}